\documentclass{amsart}
\usepackage[utf8x]{inputenc}
\usepackage{amsmath,amsthm}
\usepackage{amsfonts,amstext,amssymb, mathtools,  mathtools}

\newtheorem{proposition}{Proposition}[section]
\newtheorem{corollary}{Corollary}[section]
\newtheorem{remark}{Remark}[section]
\DeclareMathOperator{\var}{Var}
\DeclareMathOperator{\cov}{Cov}

\begin{document}
\bibliographystyle{plain}
\title{Joint asymptotic distributions of smallest and largest insurance claims}
\author{Hansj\"org Albrecher}
\address{Department of Actuarial Science, Faculty of Business and Economics, University of Lausanne, CH-1015 Lausanne, Switzerland and
Swiss Finance Institute, Switzerland}
\email{hansjoerg.albrecher@unil.ch}
\thanks{H.A. acknowledges support from the Swiss
National Science Foundation Project 200021-124635/1}

\author{Christian Y. Robert}\address{Universit\'{e} de Lyon, Universit\'{e} Lyon 1, Institut de Science Financi\`{e}re et d'Assurances, France}
\email{christian.robert@univ-lyon1.fr}

\author{Jef L.
Teugels}\address{Department of Mathematics, K.U. Leuven, Belgium}
\email{jef.teugels@wis.kuleuven.be}

\begin{abstract} Assume that claims in a portfolio of insurance contracts are described by independent and identically distributed random variables
with regularly varying tails and occur according to a near mixed Poisson process. We provide a collection of results pertaining
to the joint asymptotic Laplace transforms of the normalized sums of the smallest and
largest claims, when the length of the considered time interval tends to infinity. 
The results crucially depend on the value of the tail index of the claim distribution, as well as on the number of largest claims under consideration.
\end{abstract}
\maketitle

\section{Introduction}

When dealing with heavy-tailed insurance claims, it is a classical problem
to consider and quantify the influence of the largest among the claims on
their total sum, see e.g. Ammeter (1964) for an early reference in actuarial
literature. This topic is particularly relevant in non-proportional
reinsurance applications when a significant proportion of the sum of claims
is consumed by a small number of claims. The influence of the maximum of a
sample on the sum has in particular attracted considerable attention over
the last fifty years (see Ladoucette and Teugels \cite{LT07} for a recent
overview of existing literature on the subject). Different modes of
convergence of the ratios sum over maximum or maximum over sum have been
linked with conditions on additive domain of attractions of a stable law
(see e.g. Darling \cite{Da52}, Bobrov \cite{Bo54}, Chow and Teugels \cite%
{CT79} and Bingham and Teugels \cite{BT81}).

It is also of interest to study the joint distribution of normalized
smallest and largest claims when the number of claims over time are
described by a general counting process. This has an impact on the design of
possible reinsurance strategies and risk management in general. In this
paper we consider a homogeneous insurance portfolio, where the distribution
of the individual claims has a regularly varying tail. The number of claims
is generated by a near mixed Poisson process. For this rather general
situation we derive a number of limiting results for the joint Laplace
transforms of the smallest and largest claims, as the time $t$ tends to
infinity. These turn out to be quite explicit and crucially depend on the
rule of what is considered to be a large claim as well as on the value of
the tail index. \newline

Let $X_{1},X_{2},\ldots $ be a sequence of independent positive random
variables (representing claims) with common distribution function $F$. For $%
n\geq 1$, denote by $X_{1}^{\ast }\leq X_{2}^{\ast }\leq \ldots \leq
X_{n}^{\ast }$ the corresponding order statistics. We assume that the claim
size distribution satisfies the condition 
\begin{equation}
1-F(x)=\overline{F}(x)=x^{-\alpha }\ell(x),\qquad x>0\text{,}  \label{pareto}
\end{equation}%
where $\alpha >0$ and $\ell$ is a slowly varying function at infinity. The
tail index is defined as $\gamma =1/\alpha $ and $U(y)=F^{\leftarrow }(1-1/y)
$ is the tail quantile function of $F$. Under \eqref{pareto}, $%
U(y)=y^{1/\alpha }\ell _{1}(y)$, where $\ell_{1}$ is again a slowly varying
function. For textbook treatments of regularly varying distributions and/or
their applications in insurance modelling, see e.g. Bingham et al. \cite%
{BGT87}, Embrechts et al. \cite{EKM97}, Rolski et al. \cite{RSST99} and
Asmussen and Albrecher \cite{AA10}.

Denote the number of claims up to time $t$ by $N(t)$ with $p_{n}(t)=P(N(t)=n)
$. The probability generating function of $N(t)$ is given by 
\begin{equation*}
Q_{t}(z)=E\left\{ z^{N(t)}\right\} =\sum_{n=0}^{\infty }p_{n}(t)z^{n},
\end{equation*}%
which is defined for $|z|\leq 1$. Let 
\begin{equation*}
Q_{t}^{(r)}(z)=r!\,E\left\{ \left( 
\begin{array}{c}
N(t) \\ 
r%
\end{array}%
\right) z^{N(t)-r}\right\} 
\end{equation*}%
be its derivative of order $r$ with respect to $z$. In this paper we assume
that $N(t)$ is a near mixed Poisson (NMP) process, i.e. the claim counting
process satisfies the condition%
\begin{equation*}
\frac{N(t)}{t}\overset{D}{\rightarrow }\Theta ,\qquad t\uparrow \infty 
\end{equation*}%
for some random variable $\Theta $, where $D$ denotes convergence in
distribution. This condition implies that%
\begin{equation*}
Q_{t}\left( 1-\frac{w}{t}\right) \rightarrow E\left\{ e^{-w\Theta }\right\}
\quad \text{and}\quad \frac{1}{t^{r}}\,Q_{t}^{(r)}\left( 1-\frac{w}{t}%
\right) \rightarrow E\left\{ e^{-w\Theta }\Theta ^{r}\right\}
:=q_{r}(w),\qquad t\uparrow \infty .
\end{equation*}%
Note also that, for $\beta >0$ and $r\in \mathbb{N}$,%
\begin{equation}
\int_{0}^{\infty }w^{\beta -1}q_{r}(w)dw=\Gamma \left( \beta \right)
E\left\{ \Theta ^{r-\beta }\right\} \text{.}  \label{thiso}
\end{equation}%
If the distribution of $\Theta $ is degenerate at a single point, then $%
\left( N(t)\right) _{t\geq 0}$ has asymptotically the same behavior as a
renewal process. One particular example of a renewal process is the
homogeneous Poisson process, which is very popular in claims modelling and
plays a crucial role in both actuarial literature and practice. The general
class of NMP processes has found numerous applications in (re)insurance
modelling because of its flexibility, its success in actuarial data fitting
and its property of being more dispersed than the Poisson process (see
Grandell \cite{Gr97}). The mixing may e.g. be interpreted as claims coming
from a heterogeneity of groups of policyholders or of contract
specifications.

The aggregate claim up to time $t$ is given by 
\begin{equation*}
S(t)=\sum_{j=1}^{N(t)}X_{j},
\end{equation*}%
where it is assumed that $\left( N(t)\right) _{t\geq 0}$ is independent of
the claims $\left( X_{i}\right) _{i\geq 1}$. For $s\in \mathbb{N}$ and $%
N(t)\geq s+2$, we define the sum of the $N(t)-s-1$ smallest and the sum of
the $s$ largest claims by%
\begin{equation*}
\Sigma _{s}(t)=\sum_{j=1}^{N(t)-s-1}X_{j}^{\ast },\qquad \Lambda
_{s}(t)=\sum_{j=N(t)-s+1}^{N(t)}X_{j}^{\ast },
\end{equation*}%
so that $S(t)=\Sigma _{s}(t)+X_{N(t)-s}^{\ast }+\Lambda _{s}(t)$. Here $%
\Sigma $ refers to \textit{small} while $\Lambda $ refers to \textit{large}.\\

In this paper we study the limiting behavior of the triple $(\Lambda
_{s}(t),X_{N(t)-s}^{\ast },\Sigma _{s}(t))$ with appropriate normalisation
coefficients depending on $\gamma $, the tail index, and on $s$, the number
of terms in the sum of the largest claims. We will consider three asymptotic
cases: $s$ is fixed, $s$ tends to infinity but slower than the expected
number of claims, and $s$ tends to infinity and is asymptotically equal to a
proportion of the number of claims.\\

The paper is organized as follows. We first give the joint Laplace transform
of the triple $(\Lambda _{s}(t),X_{N(t)-s}^{\ast },\Sigma _{s}(t))$ for a
fixed $t$ in Section \ref{sec2}. Section \ref{sec3} deals with asymptotic
joint Laplace transforms in the case $0<\alpha <1$. We also discuss consequences for moments of ratios of the limiting quantities. The behavior for $\alpha=1$ 
depends on whether ${\mathbb E}[X_i]$ is finite or not. In the first case, the analysis for $\alpha>1$ applies, in the latter one has to adapt the analysis of Section \ref{sec3} exploiting the slowly varying varying function
$\int_0^xy\,dF(y)$, but we refrain from treating this very special case in detail (see e.g. \cite{AT06} for a similar adaptation in another context). Sections \ref{sec4} and %
\ref{sec5} treat the case $\alpha >1$ without and with centering,
respectively. The proofs of the results in Sections \ref{sec3}--\ref{sec5}
are given in Section \ref{sec6}. Section \ref{secc} concludes. \newline
%In Section xx we then discuss the special case of a deterministic number $n$ of claims and the limit $n\to\infty$ in more detail.

\section{Preliminaries}

\label{sec2}

In this section, we state a versatile formula that will allow us later to
derive almost all desired asymptotic properties of the joint distributions
of the triple $(\Lambda _{s}(t),X_{N(t)-s}^{\ast },\Sigma _{s}(t))$. We
consider the joint Laplace transform of $(\Lambda _{s}(t),X_{N(t)-s}^{\ast
},\Sigma _{s}(t))$ to study their joint distribution in an easy fashion. For
a fixed $t$, it is denoted by%
\begin{equation*}
\Omega _{s}(u,v,w;t)=E\left\{ \exp (-u\Lambda _{s}(t)-vX_{N(t)-s}^{\ast
}-w\Sigma _{s}(t))\right\} .
\end{equation*}%
Then the following representation holds:

\begin{proposition}
\label{Prop0}We have%
\begin{multline*}
\Omega _{s}(u,v,w;t) \\
=\sum_{n=0}^{s}p_{n}(t)\left( \int_{0}^{\infty }e^{-ux}dF(x)\right) ^{n}+%
\frac{1}{s!}\int_{0}^{\infty }\left( E[1_{\{X>y\}}e^{-uX}]\right)
^{s}e^{-vy}Q_{t}^{(s+1)}\left( E\left\{ 1_{\{X<y\}}e^{-wX}\right\} \right)
dF(y).
\end{multline*}
\end{proposition}

\noindent \textbf{Proof:} The proof is standard if we interpret $X_{r}^{\ast
}=0$ whenever $r\leq 0$. Indeed, condition on the number of claims at the
time epoch $t$ and subdivide the requested expression into three parts.%
\begin{eqnarray*}
\Omega _{s}(u,v,w;t) &=&\sum_{n=0}^{s}p_{n}(t)E\left\{ \exp \left( \left.
-u\sum_{j=1}^{n}X_{j}\right\vert N(t)=n\right) \right\} \\
&&+p_{s+1}(t)E\left\{ \exp \left( \left. -u\sum_{j=2}^{s+1}X_{j}^{\ast
}-vX_{1}^{\ast }\right\vert N(t)=s+1\right) \right\} \\
&&+\sum_{n=s+2}^{\infty }p_{n}(t)E\left\{ \exp \left( \left.
-u\sum_{j=n-s+1}^{n}X_{j}^{\ast }-vX_{n-s}^{\ast
}-w\sum_{j=1}^{n-s-1}X_{j}^{\ast }\right\vert N(t)=n\right) \right\} .
\end{eqnarray*}%
The conditional expectation in the first term on the right simplifies easily
to the form $(\int_{0}^{\infty }e^{-ux}dF(x))^{n}$. For the conditional
expectations in the second and third term, we condition additionally on the
value $y$ of the order statistic $X_{n-s}^{\ast }$; the $n-s-1$\ order
statistics $X_{1}^{\ast },X_{2}^{\ast },\ldots ,X_{n-s-1}^{\ast }$ are then
distributed independently and identically on the interval $[0,y]$ yielding
the factor $(\int_{0}^{y}e^{-wx}dF(x))^{n-s-1}$. A similar argument works
for the $s$\ order statistics $X_{n-s+1}^{\ast },X_{n-s+2}^{\ast },\ldots
,X_{n}^{\ast }$. Combinations of the two terms yields 
\begin{eqnarray*}
&&\Omega _{s}(u,v,w;t) \\
&=&\sum_{n=0}^{s}p_{n}(t)\left( \int_{0}^{\infty }e^{-ux}dF(x)\right) ^{n} \\
&&+\sum_{n=s+1}^{\infty }p_{n}(t)\frac{n!}{s!(n-s-1)!}\int_{0}^{\infty
}\left( \int_{y}^{\infty }e^{-ux}dF(x)\right) ^{s}e^{-vy}\left(
\int_{0}^{y}e^{-wx}dF(x)\right) ^{n-s-1}dF(y).
\end{eqnarray*}%
A straight-forward calculation finally shows 
\begin{eqnarray*}
&&\Omega _{s}(u,v,w;t) \\
&=&\sum_{n=0}^{s}p_{n}(t)\left( \int_{0}^{\infty }e^{-ux}dF(x)\right) ^{n}+%
\frac{1}{s!}\int_{0}^{\infty }\left( \int_{y}^{\infty }e^{-ux}dF(x)\right)
^{s}e^{-vy}Q_{t}^{(s+1)}\left( \int_{0}^{y}e^{-wx}dF(x)\right) dF(y).
\end{eqnarray*}%
\hfill $\Box $\newline

\noindent Consequently, it is possible to easily derive the expectations of
products (or ratios) of $\Lambda _{s}(t)$, $X_{N(t)-s}^{\ast }$, $\Sigma
_{s}(t)$ and $S(t)$ by differentiating (or integrating) the joint Laplace
transform. We only write down their first moment for simplicity.

\begin{corollary}
\label{Coro0}We have%
\begin{eqnarray*}
E\left\{ \Lambda _{s}(t)\right\} &=&\sum_{n=1}^{s}n\,p_{n}(t)E\left\{
X_{1}\right\} +\frac{1}{(s-1)!}\int_{0}^{\infty }\left( \overline{F}%
(y)\right) ^{s-1}\left( \int_{y}^{\infty }xdF(x)\right) Q_{t}^{(s+1)}\left(
F(y)\right) dF(y) \\
E\{X_{N(t)-s}^{\ast }\} &=&\frac{1}{s!}\int_{0}^{\infty }y\left( \overline{F}%
(y)\right) ^{s}Q_{t}^{(s+1)}\left( F(y)\right) dF(y) \\
E\left\{ \Sigma _{s}(t)\right\} &=&\frac{1}{s!}\int_{0}^{\infty }\left( 
\overline{F}(y)\right) ^{s}Q_{t}^{(s+2)}\left( F(y)\right) \left(
\int_{0}^{y}xdF(x)\right) dF(y) \\
E\left\{ S(t)\right\} &=&E\left\{ N(t)\right\} E\left\{ X_{1}\right\} .
\end{eqnarray*}
\end{corollary}

\noindent \textbf{Proof:} The individual Laplace transforms can be written
in the following form:%
\begin{eqnarray*}
E\left\{ \exp (-u\Lambda _{s}(t)\right\} &=&\sum_{n=0}^{s}p_{n}(t)\left(
\int_{0}^{\infty }e^{-ux}dF(x)\right) ^{n}+\frac{1}{s!}\int_{0}^{\infty
}\left( \int_{y}^{\infty }e^{-ux}dF(x)\right) ^{s}Q_{t}^{(s+1)}\left(
F(y)\right) dF(y) \\
E\left\{ \exp (-vX_{N(t)-s}^{\ast }\right\} &=&\Pi _{s+1}(t)+\frac{1}{s!}%
\int_{0}^{\infty }\left( \overline{F}(y)\right)
^{s}e^{-vy}Q_{t}^{(s+1)}\left( F(y\right) )dF(y) \\
E\left\{ \exp (-w\Sigma _{s}(t)\right\} &=&\Pi _{s+1}(t)+\frac{1}{s!}%
\int_{0}^{\infty }\left( \overline{F}(y)\right) ^{s}Q_{t}^{(s+1)}\left(
\int_{0}^{y}e^{-wx}dF(x)\right) dF(y) \\
E\left\{ \exp (-uS(t)\right\} &=&Q_{t}\left( \int_{0}^{\infty
}e^{-ux}dF(x)\right)
\end{eqnarray*}%
where $\Pi _{s+1}(t)=\sum_{n=0}^{s}p_{n}(t)$. By taking the first
derivative, we arrive at the respective expectations. ~\hfill $\Box $

\section{Asymptotics for the joint Laplace transforms when $0<\protect\alpha %
<1$}

\label{sec3}

Before giving the asymptotic joint Laplace transform of the sum of the
smallest and the sum of the largest claims, we first recall an important result about
convergence in distribution of order statistics and derive a
characterization of their asymptotic distribution. All proofs of this
section are deferred to Section \ref{sec6}.

It is well-known that there exists a sequence $E_{1},E_{2},...$ of
exponential random variables with unit mean such that 
\begin{equation*}
(X_{n}^{\ast },...,X_{1}^{\ast })\overset{D}{=}\left( (U\left( \Gamma
_{n+1}/\Gamma _{1}\right) ,...,U\left( \Gamma _{n+1}/\Gamma _{n}\right)
\right) 
\end{equation*}%
where $\Gamma _{k}=E_{1}+...+E_{k}$. Let $Z_{n}=\left( X_{n}^{\ast
},...,X_{1}^{\ast },0,...\right) /U(n)$. It may be shown that $Z_{n}$
converges in distribution to $Z=\left( Z_{1},Z_{2,}...\right) $ in $\mathbb{R%
}_{+}^{\mathbb{N}}$, where $Z_{k}=\Gamma _{k}^{-1/\alpha }$ (see Lemma 1 in
LePage et al. \cite{LWZ81}). For $0<\alpha <1$, the series $%
(\sum_{k=1}^{n}\Gamma _{k}^{-1/\alpha })_{n\geq 1}$ converges almost surely.
Therefore, for a fixed $s$, we deduce that, as $n\rightarrow \infty $,%
\begin{equation}
\left( \sum_{j=n-s+1}^{n}X_{j}^{\ast },X_{n-s}^{\ast
},\sum_{j=1}^{n-s-1}X_{j}^{\ast }\right) /U(n)\overset{D}{\rightarrow }%
\left( \sum_{k=1}^{s}\Gamma _{k}^{-1/\alpha },\Gamma _{s+1}^{-1/\alpha
},\sum_{k=s+2}^{\infty }\Gamma _{k}^{-1/\alpha }\right) .  \label{Lepageconv}
\end{equation}%
In particular, we derive by the Continuous Mapping Theorem that%
\begin{equation*}
\frac{\sum_{j=1}^{n-s}X_{j}^{\ast }}{X_{n-s}^{\ast }}\overset{D}{\rightarrow 
}R_{(s)}=\frac{\sum_{k=s+1}^{\infty }\Gamma _{k}^{-1/\alpha }}{\Gamma
_{s+1}^{-1/\alpha }}.
\end{equation*}%
Note that the first moment of $R_{(s)}$ (but only the first moment) may be
easily derived since%
\begin{equation}
E\left\{ R_{(s)}\right\} =1+\sum_{j=s+2}^{\infty }E\{B_{j}^{1/\alpha }\}=1+%
\frac{s+1}{\gamma -1},  \label{thist}
\end{equation}%
where $B_{j}=\sum_{i=1}^{s+1}E_{i}/\sum_{i=1}^{j}E_{i}$ has a Beta$\left(
s+1,j-1\right) $ distribution. We also recall that $F$ belongs to the
(additive) domain of attraction of a stable law with index $\alpha \in (0,1)$
if and only if 
\begin{equation*}
\lim_{n\rightarrow \infty }E\left\{ \frac{\sum_{j=1}^{n}X_{j}^{\ast }}{%
X_{n}^{\ast }}\right\} =E\left\{ R_{(0)}\right\} =1+\frac{1}{\gamma -1}=%
\frac{1}{1-\alpha }
\end{equation*}%
(see e.g. Theorem 1 in Ladoucette and Teugels \cite{LT07}).

When $\left( N(t)\right) _{t\geq 0}$ is a NMP process, we also have, as $%
t\rightarrow \infty $, 
\begin{equation}
\left( \sum_{j=N(t)-s+1}^{N(t)}X_{j}^{\ast },X_{N(t)-s}^{\ast
},\sum_{j=1}^{N(t)-s-1}X_{j}^{\ast }\right) /U(N(t))\overset{D}{\rightarrow }%
\left( \sum_{k=1}^{s}\Gamma _{k}^{-1/\alpha },\Gamma _{s+1}^{-1/\alpha
},\sum_{k=s+2}^{\infty }\Gamma _{k}^{-1/\alpha }\right)
\label{LepageconvNsto}
\end{equation}%
and 
\begin{equation*}
\frac{\sum_{j=1}^{N(t)-s}X_{j}^{\ast }}{X_{N(t)-s}^{\ast }}\overset{D}{%
\rightarrow }R_{(s)}=\frac{\sum_{k=s+1}^{\infty }\Gamma _{k}^{-1/\alpha }}{%
\Gamma _{s+1}^{-1/\alpha }}
\end{equation*}%
(see e.g. Lemma 2.5.6 in Embrechts et al. \cite{EKM97}). But note that, if
the triple $(\Lambda _{s}(t),X_{N(t)-s}^{\ast },\Sigma _{s}(t))$ is
normalized by $U(t)$ instead of $U(N(t))$ in $\left( \ref{LepageconvNsto}%
\right) $, then the asymptotic distribution will differ due to the
randomness brought in by the counting process $\left( N(t)\right) _{t\geq 0}$.\\

The following proposition gives the asymptotic Laplace transform when the
triple $(\Lambda _{s}(t),X_{N(t)-s}^{\ast },\Sigma _{s}(t))$ is normalised
by $U(t)$.

\begin{proposition}
\label{Prop1}For a fixed $s$ $\in \mathbb{N}$, as $t\rightarrow \infty $, we
have $\left( \Lambda _{s}(t)/U(t),X_{N(t)-s}^{\ast }/U(t),\Sigma
_{s}(t)/U(t)\right) \overset{D}{\rightarrow }(\Lambda _{s},\Xi _{s},\Sigma
_{s})$ where 
\begin{multline}
E\left\{ \exp (-u\Lambda _{s}-v\Xi _{s}-w\Sigma _{s})\right\}
\label{prop21e} \\
=\frac{1}{s!}\int_{0}^{\infty }\left( \frac{z}{\gamma }\int_{1}^{\infty }%
\frac{e^{-uz^{-\gamma }\eta }}{\eta ^{1+1/\gamma }}d\eta \right)
^{s}e^{-vz^{-\gamma }}q_{s+1}\left( z\left( 1+\frac{1}{\gamma }\int_{0}^{1}%
\frac{1-e^{-wz^{-\gamma }\eta }}{\eta ^{1+1/\gamma }}d\eta \right) \right)
dz.
\end{multline}%
If $\Theta =1$ a.s., this expression simplifies to%
\begin{eqnarray*}
&&E\left\{ \exp (-u\Lambda _{s}-v\Xi _{s}-w\Sigma _{s})\right\} \\
&=&\frac{1}{s!}\int_{0}^{\infty }\left( \frac{z}{\gamma }\int_{1}^{\infty }%
\frac{e^{-uz^{-\gamma }\eta }}{\eta ^{1+1/\gamma }}d\eta \right)
^{s}e^{-vz^{-\gamma }}\exp \left( -z\left( 1+\frac{1}{\gamma }\int_{0}^{1}%
\frac{1-e^{-wz^{-\gamma }\eta }}{\eta ^{1+1/\gamma }}d\eta \right) \right)
dz.
\end{eqnarray*}
\end{proposition}

We observe that $\left( N(t)\right) _{t\geq 0}$ modifies the asymptotic
Laplace transform by introducing $q_{s+1}$ into the integral $\left( \ref%
{prop21e}\right) $. However, the moments of $R_{(s)}$ do not depend on the
law of $\Theta$:

\begin{corollary}
\label{cor2} \label{Prop1moment}For $k\in \mathbb{N}^{\ast }$, we have%
\begin{equation}
E\left\{ R_{(s)}^{k}\right\} =1+\sum_{i=1}^{k}\left( 
\begin{array}{c}
k \\ 
i%
\end{array}%
\right) \sum_{j=1}^{i}\frac{(s+j)!}{s!}C_{i,j}\left( \gamma \right) .
\label{eq3b}
\end{equation}%
where%
\begin{equation}
C_{i,j}\left( \gamma \right) =\sum_{\substack{ m_{1}+\ldots +m_{i-j+1}=j  \\ %
1m_{1}+2m_{2}+\ldots +(i-j+1)m_{i-j+1}=i}}\frac{i!}{m_{1}!m_{2}!\ldots
m_{i-j+1}!}\prod_{l=1}^{i-j+1}\left( \frac{1}{l!\left( l\gamma -1\right) }%
\right) ^{m_{l}}.  \label{c-def}
\end{equation}
\end{corollary}

Note that this corollary only provides the moments of $R_{(s)}$. In order to
have moment convergence results for the ratios, it is necessary to assume
uniform integrability of $(\{\sum_{j=1}^{N(t)-s}X_{j}^{\ast
}/X_{N(t)-s}^{\ast }\}^{k})_{t\geq 0}$. It is also possible to use the
Laplace transform\ of the triple with a fixed $t$ to characterize the
moments of the ratios $\{\sum_{j=1}^{N(t)-s}X_{j}^{\ast }/X_{N(t)-s}^{\ast
}\}$ (see Corollary \ref{Coro0}), and then to follow the same approach as
proposed by Ladoucette \cite{La07} for the ratio of the random sum of
squares to the square of the random sum under the condition that $E\left\{
\Theta ^{\varepsilon }\right\} <\infty $ and $E\left\{ \Theta ^{-\varepsilon
}\right\} <\infty $ for some $\varepsilon >0$.

\begin{remark}
For $k=1$, \eqref{eq3b} reduces again to $E\left\{ R_{(s)}\right\} =1+(s+1)C_{1,1}(\gamma)=1+\frac{s+1}{\gamma-1},$ which is \eqref{thist}. Furthermore, for all $s\ge 0$
\begin{equation}\label{varian} \var\left\{ R_{(s)}\right\} =\dfrac{(s+1)\gamma ^{2}}{(\gamma
-1)^{2}(2\gamma -1)}=\dfrac{(s+1)\alpha }{(2-\alpha )(1-\alpha )^{2}}.
\end{equation}

\end{remark}

\begin{remark}
$R_{(s)}$ is the ratio of the sum $\Xi _{s}+\Sigma _{s}$ over $\Xi _{s}$. By
taking the derivative of \eqref{prop21e}, it
may be shown that, for $1<\gamma <s+1$ and $E\left\{ \Theta ^{\gamma
}\right\} <\infty $,%
\begin{equation*}
E\left\{ \Xi _{s}+\Sigma _{s}\right\} =\frac{\Gamma \left( s-\gamma
+1\right) }{(\gamma -1)\Gamma \left( s\right) }E\left\{ \Theta ^{\gamma
}\right\} .
\end{equation*}%
Therefore the mean of $\Xi _{s}+\Sigma _{s}$ will only be finite for sufficiently small $%
\gamma $. An alternative interpretation is that for given value of $\gamma$, the number $s$ of removed maximal terms in the sum has to be sufficiently large to make the mean of the remaining sum finite. 
The normalisation of the sum by $\Xi _{s}$, on the other hand, ensures the existence of the moments of the ratio $R_{(s)}$ for all values of $s$ and $\gamma>1$.
\end{remark}

\begin{remark}
It is interesting to compare Formula $\eqref{eq3b}$ with the limiting moment
of the statistic 
\begin{equation*}
T_{N(t)}=\frac{X_{1}^{2}+\cdots +X_{N(t)}^{2}}{(X_{1}+\ldots +X_{N(t)})^{2}}.
\end{equation*}%
For instance, $\lim_{t\rightarrow \infty }E\left\{ T_{N(t)}\right\}
=1-\alpha $, $\lim_{t\rightarrow \infty }\var\left\{ T_{N(t)}\right\}
=\alpha \left( 1-\alpha \right) /3$ and the limit of the $n$th moment can be
expressed as an $n$th-order polynomial in $\alpha $, see Albrecher and
Teugels \cite{AT06}, Ladoucette \cite{La07} and Albrecher et al. \cite{AST09}%
. Motivated by this similarity, let us study the link in some more detail.\
By using once again Lemma 1 in LePage et al. \cite{LWZ81}, we deduce that%
\begin{equation*}
T_{N(t)}\overset{D}{\rightarrow }T_{\infty }=\frac{\sum_{k=1}^{\infty
}\Gamma _{k}^{-2/\alpha }}{(\sum_{j=1}^{\infty }\Gamma _{k}^{-1/\alpha })^{2}%
}.
\end{equation*}%
Recall that $R_{(0)}$ is the weak limit of the ratio $%
(\sum_{j=1}^{N(t)}X_{j}^{\ast })/X_{N(t)}^{\ast }$ and $E\{R_{(0)}\}=1/(1-\alpha)$. Using \eqref{varian} and 
%$E\{R_{(0)}^{2}\}={2\gamma ^{3}}/{(\gamma -1)^{2}/(2\gamma -1)}$ and 
$E\{R_{(0)}^{2}T_{\infty }\}=2/(2-\alpha )$ (which is a
straight-forward consequence of the fact that $X_{i}^{2}$ has regularly
varying tail with tail index $2\gamma $), one then obtains a simple formula
for the covariance between $R_{(0)}^{2}$ and $T_{\infty }$: 
\begin{equation*}
\cov\left( R_{(0)}^{2},T_{\infty }\right) =-\frac{2\gamma }{1-3\gamma
+2\gamma ^{2}}.
\end{equation*}%
Determining $\var%
\{R_{(0)}^{2}\}$ by exploiting \eqref{eq3b} for $k=4$, we then arrive at the linear correlation coefficient
\begin{equation*}
\rho
(R_{(0)}^{2},T_{\infty }) =-\sqrt{\frac{3(\gamma-1)(3\gamma-1)(4\gamma-1)}{\gamma(43\gamma^2-7\gamma-6)}}.
\end{equation*}
Figure \ref{fig:corr} depicts $\rho
(R_{(0)}^{2},T_{\infty })$ as a function of $\alpha=1/\gamma$. Note that $\lim_{\gamma\to\infty}\rho
(R_{(0)}^{2},T_{\infty })=-6/\sqrt{43}$. The correlation coefficient 
allows to quantify the negative linear dependence between the two ratios
(the dependence becomes weaker when $\alpha $ increases, as the maximum term
will then typically be less dominant in the sum). % \begin{equation*}
% \rho\left( R_{(0)}^{2},T_{\infty }\right) =\frac{E\{R_{(0)}^{2}T_{\infty
% }\}-E\{R_{(0)}^{2}\}E\{T_{\infty }\}}{\sqrt{\var\{R_{(0)}^{2}\}\var\{T_{\infty }\}}%
% }=...,
% \end{equation*}%
\begin{figure}[htb]
\centering
\caption{$\rho(R_{(0)}^{2},T_{\infty })$ as a function of $\alpha$.}\label{fig:corr}
\includegraphics{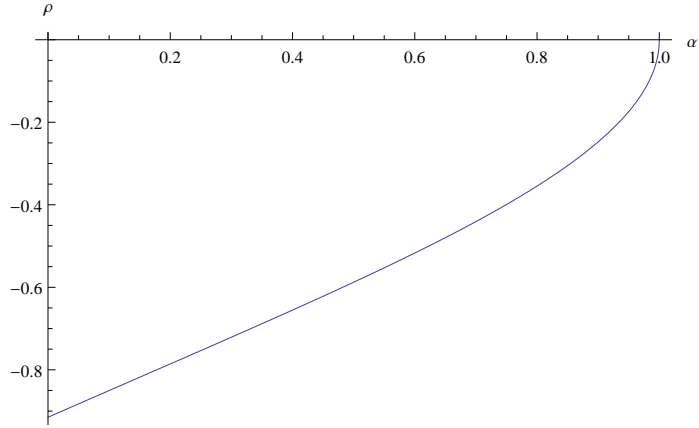}
\end{figure} 
\end{remark}

Next, let us consider the case when the number of largest terms also
increases as $t\rightarrow \infty $, but slower than the expected number of
claims. It is now necessary to change the normalisation coefficients of $%
X_{N(t)-s}^{\ast }$ and $\Sigma _{s}(t)$.

\begin{proposition}
\label{Prop2}Let $s=\left\lfloor p(t)N(t)\right\rfloor \rightarrow \infty $
for a function $p(t)$ with $p(t)\rightarrow 0$ and $tp(t)\rightarrow \infty $%
. Then $\left( \Lambda _{s}(t)/U(t),X_{N(t)-s}^{\ast }/U(p^{-1}(t)),\Sigma
_{s}(t)/(tp(t)U(p^{-1}(t)))\right) \overset{D}{\rightarrow }(\Lambda ,\Xi
,\Sigma )$ where 
\begin{equation}  \label{next}
E\left\{ \exp (-u\Lambda -v\Xi -w\Sigma )\right\} =e^{-v}q_{0}\left(
\int_{0}^{\infty }\frac{(1-e^{-uz^{-\gamma }\eta })}{\eta ^{1+1/\gamma }}%
d\eta +\frac{w}{\gamma -1}\right).
\end{equation}%
If $\Theta =1$ a.s.%
\begin{equation*}
E\left\{ \exp (-u\Lambda -v\Xi -w\Sigma )\right\} =\exp \left(
-\int_{0}^{\infty }\frac{(1-e^{-uz^{-\gamma }\eta })}{\eta ^{1+1/\gamma }}%
d\eta -v-\frac{w}{\gamma -1}\right) .
\end{equation*}%
\hfill
\end{proposition}

Several messages may be derived from (\ref{next}). First note that the
asymptotic distribution of $X_{N(t)-s}^{\ast }$ is degenerated for $%
s=\left\lfloor p(t)N(t)\right\rfloor$, since $X_{N(t)-s}^{\ast }/U(p^{-1}(t))%
\overset{D}{\rightarrow }1$ as $t\rightarrow \infty $. Second, the
asymptotic distribution of the sum of the smallest claims is the
distribution of $\Theta $ up to a scaling factor, since $\Sigma
_{s}(t)/(tp(t)U(p^{-1}(t)))\overset{D}{\rightarrow }\Theta /(\gamma -1)$ as $%
t\rightarrow \infty $.\newline

Finally, for a fixed proportion of maximum terms, it is also necessary to
change the normalisation coefficients of $X_{N(t)-s}^{\ast }$ and $\Sigma
_{s}(t)$. We have

\begin{proposition}
\label{Prop3}Let $s=\left\lfloor pN(t)\right\rfloor $ for a fixed $0<p<1$.
Then $\left( \Lambda _{s}(t)/U(t),X_{N(t)-s}^{\ast },\Sigma _{s}(t)/t\right) 
\overset{D}{\rightarrow }(\Lambda _{p},\Xi _{p},\Sigma _{p})$ where%
\begin{equation*}
E\left\{ \exp (-u\Lambda _{p}-v\Xi _{p}-w\Sigma _{p})\right\}
=e^{-vx_{p}}q_{0}\left( u^{\alpha }\frac{\Gamma (1-\alpha )}{1-p}+wE\left\{
X|X\leq x_{p}\right\} \right) 
\end{equation*}%
and $x_{p}=F^{-1}(p)$. If $\Theta =1$ a.s.,%
\begin{equation*}
E\left\{ \exp (-u\Lambda _{p}-v\Xi _{p}-w\Sigma _{p})\right\} =\exp \left(
-u^{\alpha }\frac{\Gamma (1-\alpha )}{1-p}-vx_{p}-wE\left\{ X|X\leq
x_{p}\right\} \right) .
\end{equation*}
\end{proposition}

As expected, $X_{N(t)-s}^{\ast }\overset{D}{\rightarrow }x_{p}$ and $\Sigma
_{s}(t)/t\overset{D}{\rightarrow }\Theta E\left\{ X|X\leq x_{p}\right\} $ as 
$t\rightarrow \infty $. If $\Theta =1$ a.s. and $\alpha =1/2$, then $\Lambda
_{p}$ has an inverse Gamma distribution with shape parameter equal to $1/2$.

\section{Asymptotics for the joint Laplace transforms when $\protect\alpha %
>1 $}

\label{sec4}

In this section, we assume that $\alpha >1$ and hence the expectation of the
claim distribution is finite. We let $\mu =E\left\{ X_{1}\right\} $. The
normalisation coefficient of the sum of the smallest claims, $\Sigma _{s}(t)$%
, will therefore be $t^{-1}$ as it is the case for $S(t)$ for the Law of
Large Numbers. In Section \ref{sec5}, we will then consider the sum of the
smallest centered claims with another normalisation coefficient.

Again, consider fixed {$s$ $\in \mathbb{N}$} first. The normalisation
coefficients of $\Lambda _{s}(t)$ and $X_{N(t)-s}^{\ast }$ are the same as
for the case $0<\alpha <1$, but the normalisation coefficient of $\Sigma _{s}
$ is now $t^{-1}$.

\begin{proposition}
\label{Prop4} For fixed $s$ $\in \mathbb{N}$, we have $\left( \Lambda
_{s}(t)/U(t),X_{N(t)-s}^{\ast }/U(t),\Sigma _{s}(t)/t\right) \overset{D}{%
\rightarrow }(\Lambda _{s},\Xi _{s},\Sigma _{s})$ where 
\begin{equation*}
E\left\{ \exp (-u\Lambda _{s}-v\Xi _{s}-w\Sigma _{s})\right\} =\frac{1}{s!}%
\int_{0}^{\infty }\left( \frac{z}{\gamma }\int_{1}^{\infty }\frac{%
e^{-uz^{-\gamma }\eta}}{\eta^{1+1/\gamma }}d\eta\right) ^{s}e^{-vz^{-\gamma
}}q_{s+1}\left( z+w\mu \right) dz.
\end{equation*}%
If $\Theta =1$ a.s.,%
\begin{equation*}
E\left\{ \exp (-u\Lambda _{s}-v\Xi _{s}-w\Sigma _{s})\right\} =e^{-w\mu }%
\frac{1}{s!}\int_{0}^{\infty }\left( \frac{z}{\gamma }\int_{1}^{\infty }%
\frac{e^{-uz^{-\gamma }\eta}}{\eta^{1+1/\gamma }}d\eta\right)
^{s}e^{-vz^{-\gamma }-z}dz.
\end{equation*}
\end{proposition}

\begin{corollary}
\label{Prop4moment}We have%
\begin{equation*}
E\left\{ \frac{\Sigma _{s}}{\Xi _{s}}\right\} =\mu \frac{\Gamma (s-\gamma +1)%
}{s!}E\left\{ \Theta ^{1+\gamma }\right\}
\end{equation*}%
and%
\begin{equation*}
E\left\{ \frac{\Xi _{0}}{\Lambda _{s}+\Xi _{s}+\Sigma _{s}}\right\} =1-\mu
\int_{0}^{\infty }\int_{0}^{\infty }e^{-uz^{-\gamma }}q_{2}\left( z+u\mu
\right) dudz.
\end{equation*}
\end{corollary}

We first note that 
\begin{equation*}
E\left\{ \exp (-w\Sigma _{s})\right\} =\frac{1}{s!}\int_{0}^{\infty
}z^{s}q_{s+1}\left( z+w\mu \right) dz=E\left\{ \frac{1}{s!}\int_{0}^{\infty
}z^{s}(e^{-\left( z+w\mu \right) \Theta }\Theta ^{s+1})dz\right\} =E\left\{
e^{-w\mu \Theta }\right\} 
\end{equation*}%
and therefore $\Sigma _{s}(t)/t\overset{D}{\rightarrow }\mu \Theta $ as $%
t\rightarrow \infty $ for any fixed {$s$ $\in \mathbb{N}$. }The influence of
the largest claims on the sum becomes less and less important as $t$ is
large and is asymptotically negligible. This is very different from the case 
$0<\alpha <1$. In Theorem 1 in Downey and Wright \cite{DW07}, it is moreover
shown that, as $n\rightarrow \infty $, 
\begin{equation*}
E\left\{ \frac{X_{n}^{\ast }}{\sum_{j=1}^{n}X_{j}^{\ast }}\right\} =\frac{%
E\left\{ X_{n}^{\ast }\right\} }{E\{\sum_{j=1}^{n}X_{j}^{\ast }\}}\left(
1+o(1)\right) .
\end{equation*}%
This result is no more true in our framework when $\Theta $ is not
degenerate at $1$. Assume that $E\left\{ \Theta ^{\gamma }\right\} <\infty $%
. Using \eqref{thiso} and under a uniform integrability condition, one has%
\begin{equation*}
\lim_{t\rightarrow \infty }E\left\{ \frac{X_{N(t)}^{\ast }}{%
\sum_{j=1}^{N(t)}X_{j}^{\ast }}\right\} \frac{t}{U(t)}\neq \frac{%
\lim_{t\rightarrow \infty }E\left\{ X_{N(t)}^{\ast }\right\} /U(t)}{%
\lim_{t\rightarrow \infty }E\left\{ \sum_{j=1}^{N(t)}X_{j}^{\ast }\right\} /t%
}=\frac{\Gamma (1-\gamma )E\left\{ \Theta ^{\gamma }\right\} }{\mu E\left\{
\Theta \right\} }.
\end{equation*}

Next, we consider the case with varying number of maximum terms. The
normalisation coefficients of $\Lambda _{s}(t)$ and $X_{N(t)-s}^{\ast }$ now
differ.

\begin{proposition}
\label{Prop5}Let {$s=\left\lfloor p(t)N(t)\right\rfloor \rightarrow \infty $
and $p(t)\rightarrow 0$, i.e. $tp(t)\rightarrow \infty $}. Then $$\left(
\Lambda _{s}(t)/(tp(t)U(p^{-1}(t))),X_{N(t)-s}^{\ast }/U(p^{-1}(t)),\Sigma
_{s}(t)/t\right) \overset{D}{\rightarrow }(\Lambda ,\Xi ,\Sigma ),$$ where%
\begin{equation*}
E\left\{ \exp (-u\Lambda -v\Xi -w\Sigma )\right\} =e^{-v}q_{0}\left( \frac{u%
}{1-\gamma }+w\mu \right) .
\end{equation*}%
If $\Theta =1$ a.s.,%
\begin{equation*}
E\left\{ \exp (-u\Lambda -v\Xi -w\Sigma )\right\} =e^{-u/(1-\gamma
)}e^{-v}e^{-w\mu }.
\end{equation*}
\end{proposition}

As for the case $0<\alpha <1$, $X_{N(t)-s}^{\ast }/U(p^{-1}(t))\overset{P}{%
\rightarrow }1$ as $t\rightarrow \infty $. Moreover the asymptotic
distribution of the sum of the largest claim is the distribution of $\Theta $
up to a scaling factor since $\Lambda _{s}(t)/(tp(t)U(p^{-1}(t)))\overset{D}{%
\rightarrow }\Theta /(1-\gamma )$ as $t\rightarrow \infty $. Finally note
that $\Sigma _{s}(t)/t\overset{D}{\rightarrow }\mu \Theta $ as $t\rightarrow
\infty $ as for the case when $s$ was fixed.\newline

Finally we fix $p$. Only the normalisation coefficient of $\Lambda _{s}(t)$
and its asymptotic distribution differ from the case $0<\alpha <1$.

\begin{proposition}
\label{Prop6}Let {$s=\left\lfloor pN(t)\right\rfloor $ and $0<p<1$}. Then $%
\left( \Lambda _{s}(t)/t,X_{N(t)-s}^{\ast },\Sigma _{s}(t)/t\right) \overset{%
D}{\rightarrow }(\Lambda _{p},\Xi _{p},\Sigma _{p})$ where%
\begin{equation*}
E\left\{ \exp (-u\Lambda _{p}-v\Xi _{p}-w\Sigma _{p})\right\}
=e^{-vx_{p}}q_{0}\left( uE\left\{ X|X>x_{p}\right\} +wE\left\{ X|X\leq
x_{p}\right\} \right)
\end{equation*}%
and $x_{p}=F^{-1}(p)$. If $\Theta =1$ a.s.,%
\begin{equation*}
E\left\{ \exp (-u\Lambda _{p}-v\Xi _{p}-w\Sigma _{p})\right\} =e^{-uE\left\{
X|X>x_{p}\right\} }e^{-vx_{p}}e^{-wE\left\{ X|X\leq x_{p}\right\} }.
\end{equation*}
\end{proposition}

We note that the normalisation of $\Lambda _{s}(t)$ is the same as for $%
\Sigma _{s}(t)$ and that $\Lambda _{s}(t)/t\overset{D}{\rightarrow }\Theta
E\left\{ X|X>x_{p}\right\} $ as $t\rightarrow \infty $.

\section{Asymptotics for the joint Laplace transform for $\protect\alpha >1$
with centered claims when $s$ is fixed}

\label{sec5}

In this section, we consider the sum of the smallest centered claims:%
\begin{equation*}
\Sigma _{s}^{(\mu )}(t)=\sum_{j=1}^{N(t)-s-1}\left( X_{j}^{\ast }-\mu
\right) .
\end{equation*}%
instead of the sum of the smallest claims $\Sigma _{s}(t)$. Like for the
Central Limit Theorem, we have to consider two subcases: $1<\alpha <2$ and {$%
\alpha >2$}${.}$

For the subcase $1<\alpha <2$, the normalisation coefficient of $\Sigma
_{s}^{(\mu )}(t)$ is now $U^{-1}(t)$.

\begin{proposition}
\label{Prop7}For fixed $s$ $\in \mathbb{N}$ and $1<\alpha <2$, we have 
\begin{equation*}
\left( \Lambda _{s}(t)/U(t),X_{N(t)-s}^{\ast }/U(t),\Sigma _{s}^{(\mu
)}(t)/U(t)\right) \overset{D}{\rightarrow }(\Lambda _{s},\Xi _{s},\Sigma
_{s}^{(\mu )}),
\end{equation*}%
where 
\begin{eqnarray*}
&&E\left\{ \exp (-u\Lambda _{s}-v\Xi _{s}-w\Sigma _{s}^{(\mu )})\right\}  \\
&=&\frac{1}{s!}\int_{0}^{\infty }\left( \frac{z}{\gamma }\int_{1}^{\infty }%
\frac{e^{-uz^{-\gamma }\eta }}{\eta ^{1+1/\gamma }}d\eta \right)
^{s}e^{-vz^{-\gamma }}q_{s+1}\left( z\left( 1+\frac{1}{\gamma }\int_{0}^{1}%
\frac{1-wz^{-\gamma }\eta -e^{-wz^{-\gamma }\eta }}{\eta ^{1+1/\gamma }}%
d\eta -\frac{z^{-\gamma }}{1-\gamma }w\right) \right) dz.
\end{eqnarray*}%
If $\Theta =1$ a.s.,%
\begin{eqnarray*}
&&E\left\{ \exp (-u\Lambda _{s}-v\Xi _{s}-w\Sigma _{s}^{(\mu )})\right\}  \\
&=&\frac{1}{s!}\int_{0}^{\infty }\left( \frac{z}{\gamma }\int_{1}^{\infty }%
\frac{e^{-uz^{-\gamma }\eta }}{\eta ^{1+1/\gamma }}d\eta \right)
^{s}e^{-vz^{-\gamma }}\exp \left( -z\left( 1+\frac{1}{\gamma }\int_{0}^{1}%
\frac{1-wz^{-\gamma }\eta -e^{-wz^{-\gamma }\eta }}{\eta ^{1+1/\gamma }}%
d\eta -\frac{z^{-\gamma }}{1-\gamma }w\right) \right) dz.
\end{eqnarray*}
\end{proposition}

If $s=0$, then 
\begin{equation*}
E\left\{ \exp (-v\Xi _{0}-w\Sigma _{0}^{(\mu )})\right\} =\int_{0}^{\infty
s}e^{-vz^{-\gamma }}\exp \left( -z\left( 1+\frac{1}{\gamma }\int_{0}^{1}%
\frac{1-wz^{-\gamma }\eta -e^{-wz^{-\gamma }\eta }}{\eta ^{1+1/\gamma }}%
d\eta -\frac{z^{-\gamma }}{1-\gamma }w\right) \right) dz
\end{equation*}%
and we see that $\Xi _{0}$ and $\Sigma _{0}^{(\mu )}$ are not independent.

\begin{corollary}
\label{Prop7moment}We have%
\begin{equation*}
E\left\{ 1+\frac{\Sigma _{s}^{(\mu )}}{\Xi _{s}}\right\} =1+\frac{s+1}{%
\gamma -1}.
\end{equation*}
\end{corollary}

This result is to compare with the one obtained by Bingham and Teugels \cite%
{BT81} for $s=0$ (see also Ladoucette and Teugels \cite{LT07}).

For the subcase {$\alpha >2$}, let $\sigma ^{2}=\var\left\{ X_{1}\right\} $.
The normalisation coefficient of $\Sigma _{s}^{(\mu )}(t)$ becomes $t^{-1/2}$%
.

\begin{proposition}
\label{Prop8}For {$s$ $\in \mathbb{N}$ fixed and $\alpha >2$}, we have $%
\left( \Lambda _{s}(t)/U(t),X_{N(t)-s}^{\ast }/U(t),\Sigma _{s}^{(\mu
)}(t)/t^{1/2}\right) \overset{D}{\rightarrow }(\Lambda _{s},\Xi _{s},\Sigma
_{s}^{(\mu )})$ where 
\begin{equation*}
E\left\{ \exp (-u\Lambda _{s}-v\Xi _{s}-w\Sigma _{s}^{(\mu )})\right\} =%
\frac{1}{s!}\int_{0}^{\infty }\left( \frac{z}{\gamma }\int_{1}^{\infty }%
\frac{e^{-uz^{-\gamma }\eta}}{\eta^{1+1/\gamma }}d\eta\right)
^{s}e^{-vz^{-\gamma }}q_{s+1}\left( z-\frac{1}{2}w^{2}\sigma ^{2}\right) dz
\end{equation*}%
If $\Theta =1$ a.s.,%
\begin{equation*}
E\left\{ \exp (-u\Lambda _{s}-v\Xi _{s}-w\Sigma _{s}^{(\mu )})\right\} =%
\frac{1}{s!}\,\exp \left( \frac{1}{2}w^{2}\sigma ^{2}\right)\int_{0}^{\infty
}\left( \frac{z}{\gamma }\int_{1}^{\infty }\frac{e^{-uz^{-\gamma }\eta}}{%
\eta^{1+1/\gamma }}d\eta\right) ^{s}e^{-vz^{-\gamma }-z}dz .
\end{equation*}
\end{proposition}

If $s=0$ and $\Theta =1$ a.s., we note that the maximum, $\Xi _{0}$, and the
centered sum, $\Sigma _{0}^{(\mu )}$, are independent. If $s>0$ and $\Theta
=1$ a.s., $\left( \Lambda _{s},\Xi _{s}\right) $ is independent of $\Sigma
_{s}^{(\mu )}$.

\section{Proofs}

\label{sec6} % 

\noindent\textbf{Proof of Proposition \ref{Prop1}:} In formula %
\eqref{prop21e}, we first use the substitution $\overline{F}(y)=z/t$, i.e. $%
y=U(t/z)$:\ 
\begin{eqnarray*}
&&\Omega _{s}(u/U(t),v/U(t),w/U(t)) \\
&=&\sum_{n=0}^{s}p_{n}(t)\left( \int_{0}^{\infty }e^{-ux/U(t)}dF(x)\right)
^{n} \\
&&+\frac{1}{s!}\int_{0}^{t}\left( \int_{U(t/z)}^{\infty
}e^{-ux/U(t)}dF(x)\right) ^{s}e^{-vU(t/z)/U(t)}Q_{t}^{(s+1)}\left(
\int_{0}^{U(t/z)}e^{-wx/U(t)}dF(x)\right) \frac{dz}{t} 
\end{eqnarray*}
\begin{eqnarray*}
&=&\sum_{n=0}^{s}p_{n}(t)\left( \int_{0}^{\infty }e^{-ux/U(t)}dF(x)\right)
^{n} \\
&&+\frac{1}{s!}\int_{0}^{t}\left( t\int_{U(t/z)}^{\infty
}e^{-ux/U(t)}dF(x)\right) ^{s}e^{-vU(t/z)/U(t)} \\
&&\times \frac{1}{t^{s+1}}Q_{t}^{(s+1)}\left( 1-\frac{1}{t}\left(
t-t\int_{0}^{U(t/z)}e^{-wx/U(t)}dF(x)\right) \right) dz.
\end{eqnarray*}%
Next, the substitution $\overline{F}(x)=\rho z/t$, i.e. $x=U(t/(z\rho ))$
leads to%
\begin{equation*}
t\int_{U(t/z)}^{\infty }e^{-ux/U(t)}dF(x)=z\int_{0}^{1}e^{-uU(t/(z\rho
))/U(t)}d\rho \rightarrow z\int_{0}^{1}e^{-u(z\rho )^{-\gamma }}d\rho =\frac{%
z}{\gamma }\int_{1}^{\infty }\frac{e^{-uz^{-\gamma }\eta}}{\eta^{1+1/\gamma }%
}\,d\eta
\end{equation*}%
as ${t\to\infty}$ and also%
\begin{multline*}
\left( t-t\int_{0}^{U(t/z)}e^{-wx/U(t)}dF(x)\right)
=t(1-F(U(t/z)))+t\int_{0}^{U(t/z)}(1-e^{-wx/U(t)})dF(x) \\
=z+z\int_{1}^{\infty }(1-e^{-wU(t/(z\rho ))/U(t)})d\rho \\
\rightarrow z\left( 1+\int_{1}^{\infty }(1-e^{-w(z\rho )^{-\gamma }})d\rho
\right) =z\left( 1+\frac{1}{\gamma }\int_{0}^{1}\frac{1-e^{-wz^{-\gamma
}\eta}}{\eta^{1+1/\gamma }}d\eta\right).
\end{multline*}%
Note that the integral is well defined since $\gamma >1$. Moreover $%
e^{-vU(t/z)/U(t)}\rightarrow e^{-vz^{-\gamma }}$ and%
\begin{equation*}
p_{n}(t)\left( \int_{0}^{\infty }e^{-ux/U(t)}dF(x)\right) ^{n}\leq
p_{n}(t)\rightarrow 0\quad\text{as}\;t\to\infty.
\end{equation*}%
\hfill $\Box$\newline

\textbf{Proof of Corollary \ref{Prop1moment}}: % By Proposition \ref{Prop1}
From Proposition \ref{Prop1} we have 
\begin{equation*}
E\left\{ \exp (-(u+v)\Xi _{s}-u\Sigma _{s})\right\} =\frac{1}{s!}%
\int_{0}^{\infty }z^{s}e^{-uz^{-\gamma }}e^{-vz^{-\gamma }}q_{s+1}\left(
z\left( 1+\frac{1}{\gamma }\int_{0}^{1}\frac{1-e^{-uz^{-\gamma }\eta }}{\eta
^{1+1/\gamma }}d\eta \right) \right) dz.
\end{equation*}%
Hence%
\begin{eqnarray*}
\left. \frac{\partial }{\partial u}E\left\{ \exp (-(u+v)\Xi _{s}-u\Sigma
_{s})\right\} \right\vert _{u=0} &=&-\frac{1}{s!}\int_{0}^{\infty
}z^{s}z^{-\gamma }e^{-vz^{-\gamma }}q_{s+1}\left( z\right) dz \\
&&-\frac{1}{s!\left( \gamma -1\right) }\int_{0}^{\infty }z^{s+1}z^{-\gamma
}e^{-vz^{-\gamma }}q_{s+2}\left( z\right) dz.
\end{eqnarray*}%
This gives indeed, using (\ref{thiso}), 
\begin{multline*}
E\left\{ \frac{\Xi _{s}+\Sigma _{s}}{\Xi _{s}}\right\} =-\int_{0}^{\infty
}\left. \frac{\partial }{\partial u}E\left\{ \exp (-(u+v)\Xi _{s}-u\Sigma
_{s})\right\} \right\vert _{u=0}dv \\
=\frac{1}{s!}\int_{0}^{\infty }z^{s}q_{s+1}\left( z\right) dz+\frac{1}{%
s!\left( \gamma -1\right) }\int_{0}^{\infty }z^{s+1}q_{s+2}\left( z\right)
dz=1+\frac{s+1}{\gamma -1},
\end{multline*}%
which extends (\ref{thist}) to the case of NMP processes. Next, we focus on %
\eqref{eq3b} for general $k$. We first consider the case $s=0$. We have%
\begin{equation*}
E\left\{ \left( 1+\frac{\Sigma _{0}}{\Xi _{0}}\right) ^{k}\right\}
=\sum_{i=0}^{k}\left( 
\begin{array}{c}
k \\ 
i%
\end{array}%
\right) E\left\{ \left( \frac{\Sigma _{0}}{\Xi _{0}}\right) ^{i}\right\} .
\end{equation*}%
Let%
\begin{equation*}
\theta (z,w)=z\left( 1+\frac{1}{\gamma }\int_{0}^{1}\frac{1-e^{-wz^{-\gamma
}\eta }}{\eta ^{1+1/\gamma }}d\eta \right) .
\end{equation*}%
By Proposition \ref{Prop1}

\begin{equation}
E\left\{ \exp (-v\Xi _{0}-w\Sigma _{0})\right\} =\int_{0}^{\infty
}e^{-vz^{-\gamma }}q_{1}\left( \theta (z,w)\right) dz  \label{p0}
\end{equation}%
and clearly%
\begin{equation*}
E\left\{ \left( \frac{\Sigma _{0}}{\Xi _{0}}\right) ^{i}\right\} =\frac{%
(-1)^{i}}{\Gamma (i)}\int_{0}^{\infty }v^{i-1}\left. \frac{\partial ^{i}}{%
\partial w^{i}}E(\exp (-v\Xi _{0}-w\Sigma _{0}))\right\vert _{w=0}dv.
\end{equation*}

Note that%
\begin{eqnarray*}
\theta ^{\left( 1\right) }(z,w):= &&\frac{\partial }{\partial w}\theta (z,w)=%
\frac{z^{-\gamma +1}}{\gamma }\int_{0}^{1}\eta ^{-1/\gamma }e^{-wz^{-\gamma
}\eta }d\eta  \\
\theta ^{\left( n\right) }(z,w):= &&\frac{\partial ^{n}}{\partial w^{n}}%
\theta (z,w)=\left( -1\right) ^{n+1}\frac{z^{-n\gamma +1}}{\gamma }%
\int_{0}^{1}\eta ^{-1/\gamma +(n-1)}e^{-wz^{-\gamma }\eta }d\eta ,
\end{eqnarray*}%
so 
\begin{equation*}
\theta ^{\left( n\right) }(z,0)=\left( -1\right) ^{n+1}z^{-n\gamma +1}\frac{1%
}{n\gamma -1}.
\end{equation*}%
By de Faa di Bruno's formula%
\begin{equation*}
\frac{\partial ^{n}}{\partial w^{n}}q_{1}\left( \theta (z,w)\right)
=\sum_{k=0}^{n}q_{1}^{(k)}\left( \theta (z,w)\right) B_{n,k}\left( \theta
^{\left( 1\right) }(z,w),\ldots ,\theta ^{\left( n-k+1\right) }(z,w)\right) 
\end{equation*}%
where%
\begin{equation*}
B_{n,k}\left( x_{1},\ldots ,x_{n-k+1}\right) =\sum_{\substack{ m_{1}+\ldots
+m_{n-k+1}=k \\ 1m_{1}+2m_{2}+\ldots +(n-k+1)m_{n-k+1}=n}}\frac{n!}{%
m_{1}!m_{2}!\ldots m_{n-k+1}!}\prod_{j=1}^{n-k+1}\left( \frac{x_{j}}{j!}%
\right) ^{m_{j}}.
\end{equation*}%
Therefore%
\begin{equation*}
\left. \frac{\partial ^{n}}{\partial w^{n}}q_{1}\left( \theta (z,w)\right)
\right\vert _{w=0}=\sum_{k=1}^{n}\left( -1\right) ^{k}q_{k+1}\left( z\right)
B_{n,k}\left( \frac{z^{-\gamma +1}}{\gamma -1},\ldots ,\left( -1\right)
^{n-k}\frac{z^{-(n-k+1)\gamma +1}}{(n-k+1)\gamma -1}\right) .
\end{equation*}%
Subsequently, 
\begin{eqnarray*}
&&B_{n,k}\left( \frac{z^{-\gamma +1}}{\gamma -1},\ldots ,\left( -1\right)
^{n-k}\frac{z^{-(n-k+1)\gamma +1}}{(n-k+1)\gamma -1}\right)  \\
&=&\sum_{\substack{ m_{1}+\ldots +m_{n-k+1}=k \\ 1m_{1}+2m_{2}+\ldots
+(n-k+1)m_{n-k+1}=n}}\frac{n!}{m_{1}!m_{2}!\ldots m_{n-k+1}!}%
\prod_{j=1}^{n-k+1}\left( \frac{\left( -1\right) ^{j+1}}{j!}z^{-j\gamma +1}%
\frac{1}{j\gamma -1}\right) ^{m_{j}} \\
&=&z^{-n\gamma +k}\left( -1\right) ^{n+k}\sum_{\substack{ m_{1}+\ldots
+m_{n-k+1}=k \\ 1m_{1}+2m_{2}+\ldots +(n-k+1)m_{n-k+1}=n}}\frac{n!}{%
m_{1}!m_{2}!\ldots m_{n-k+1}!}\prod_{j=1}^{n-k+1}\left( \frac{1}{j!\left(
j\gamma -1\right) }\right) ^{m_{j}} \\
&=&z^{-n\gamma +k}\left( -1\right) ^{n+k}C_{n,k}\left( \gamma \right) 
\end{eqnarray*}%
with definition \eqref{c-def}. This gives%
\begin{equation*}
\left. \frac{\partial ^{n}}{\partial w^{n}}q_{1}\left( \theta (z,w)\right)
\right\vert _{w=0}=z^{-n\gamma +k}\left( -1\right)
^{n}\sum_{k=1}^{n}q_{k+1}\left( z\right) C_{n,k}\left( \gamma \right) .
\end{equation*}%
and%
\begin{equation*}
\left. \frac{\partial ^{i}}{\partial w^{i}}E(\exp (-v\Xi _{0}-w\Sigma
_{0}))\right\vert _{w=0}=\int_{0}^{\infty }e^{-vz^{-\gamma }}\left(
z^{-i\gamma +k}\left( -1\right) ^{i}\sum_{k=1}^{i}q_{k+1}\left( z\right)
C_{i,k}\left( \gamma \right) \right) dz
\end{equation*}%
so that%
\begin{eqnarray*}
E\left\{ \left( \frac{\Sigma _{0}}{\Xi _{0}}\right) ^{i}\right\}  &=&\frac{1%
}{\Gamma (i)}\int_{0}^{\infty }v^{i-1}\left[ \int_{0}^{\infty
}e^{-vz^{-\gamma }}\left( z^{-i\gamma +k}\sum_{k=1}^{i}q_{k+1}\left(
z\right) C_{i,k}\left( \gamma \right) \right) dz\right] dv \\
&=&\frac{1}{\Gamma (i)}\sum_{k=1}^{i}C_{i,k}\left( \gamma \right)
\int_{0}^{\infty }q_{k+1}\left( z\right) z^{-i\gamma +k}\left[
\int_{0}^{\infty }v^{i-1}e^{-vz^{-\gamma }}dv\right] dz \\
&=&\frac{1}{\Gamma (i)}\sum_{k=1}^{i}C_{i,k}\left( \gamma \right)
\int_{0}^{\infty }q_{k+1}\left( z\right) z^{-i\gamma +k}\frac{\Gamma (i)}{%
\left( z^{-\gamma }\right) ^{i}}dz \\
&=&\sum_{k=1}^{i}C_{i,k}\left( \gamma \right) \int_{0}^{\infty
}q_{k+1}\left( z\right) z^{k}dz \\
&=&\sum_{k=1}^{i}k!C_{i,k}\left( \gamma \right) ,
\end{eqnarray*}%
cf. \eqref{thiso}, and the result follows.

For the case $s>0$, we proceed in an analogous way. Equation \eqref{p0}
becomes 
\begin{equation*}
E\left\{ \exp (-v\Xi _{s}-w\Sigma _{s})\right\} =\int_{0}^{\infty
}z^{s}e^{-vz^{-\gamma }}q_{s+1}\left( \theta (z,w)\right) dz.
\end{equation*}%
Then $\Sigma _{0}/\Xi _{0}$ is replaced by $\Sigma _{s}/\Xi _{s}$, $%
q_{1}\left( z\right) $ by $q_{s+1}\left( z\right) $, $q_{1}^{(k)}$ by $%
q_{s+1}^{(k)}$ and, by following the same path as for $s=0$, we get%
\begin{equation*}
E\left\{ \left( \frac{\Sigma _{s}}{\Xi _{s}}\right) ^{i}\right\}
=\sum_{k=1}^{i}C_{i,k}\left( \gamma \right) \int_{0}^{\infty
}q_{s+k+1}\left( z\right) z^{s+k}dz=\sum_{k=1}^{i}\frac{(s+k)!}{s!}%
C_{i,k}\left( \gamma \right) .
\end{equation*}
\hfill $\Box $\newline

\textbf{Proof of Proposition \ref{Prop2}:} The proof is similar to the
previous one, so we just highlight the differences here: Conditioning on $%
N(t)=n$ we have 
\begin{eqnarray*}
&&E_{N(t)=n}\left\{ \exp (-u\Lambda _{s}(t)/U(t)-vX_{n(1-p(t))}^{\ast
}/U(p^{-1}(t)))-w\Sigma _{s}(t))/(tp(t)U(p^{-1}(t)))\right\} \\
&=&\frac{n!}{(np(t))!(n(1-p(t)))!}\int_{0}^{\infty }\left( \int_{y}^{\infty
}e^{-ux/U(t)}dF(x)\right) ^{np(t)}e^{-vy/U(p^{-1}(t))} \\
&&\times \left( \int_{0}^{y}e^{-wx/(tp(t)U(p^{-1}(t)))}dF(x)\right)
^{n(1-p(t))-1}dF(y).
\end{eqnarray*}%
We first replace $F$ by the substitution $\overline{F}(y)=p(t)z$, i.e. $%
y=U(1/(p(t)z))$ :\textbf{\ }\ 
\begin{eqnarray*}
&&\left. \Omega
_{s}(u/U(t),v/U(p^{-1}(t)),w/(tp(t)U(p^{-1}(t)));t)\right\vert _{N(t)=n} \\
&=&\frac{n!}{(np(t))!(n(1-p(t)))!}\int_{0}^{\infty }\left(
\int_{U(1/(p(t)z))}^{\infty }e^{-ux/(tp(t)U(p^{-1}(t))}dF(x)\right)
^{np(t)}e^{-vU(1/(p(t)z))/U(p^{-1}(t))} \\
&&\times \left(
\int_{0}^{U(1/(p(t)z))}e^{-wx/(tp(t)U(p^{-1}(t)))}dF(x)\right)
^{n(1-p(t))-1}p(t)dz.
\end{eqnarray*}%
The factor involving $v$ converges to 
\begin{equation*}
e^{-vU(1/(p(t)z))/U(p^{-1}(t))}\rightarrow e^{-vz^{-\gamma }}.
\end{equation*}%
The factor containing $w$ behaves as 
\begin{eqnarray*}
\int_{0}^{U(1/(p(t)z))}e^{-wx/(tp(t)U(p^{-1}(t)))}dF(x)
&=&1-p(t)z-\int_{0}^{U(1/(p(t)z))}(1-e^{-wx/(tp(t)U(p^{-1}(t)))}dF(x) \\
&=&1-p(t)z-w\int_{0}^{U(1/(p(t)z))}\frac{x}{(tp(t)U(p^{-1}(t)))}dF(x)+... \\
&=&1-p(t)z-w\frac{r\left( U(1/(p(t)z))\right) }{(tp(t)U(p^{-1}(t)))}+... \\
&=&1-p(t)z-\frac{w}{t}\frac{z^{1-\gamma }}{\gamma -1}+...
\end{eqnarray*}%
and hence for the power%
\begin{multline*}
\left( \int_{0}^{U(1/(p(t)z))}e^{-wx/(tp(t)U(p^{-1}(t)))}dF(x)\right)
^{n(1-p(t))-1} \\
=\exp \left( n(1-p(t))\ln \left( 1-p(t)z-\frac{w}{t}\frac{z^{1-\gamma }}{%
\gamma -1}+...\right) \right) \\
=\exp \left( np(t)z-\frac{n}{t}w\frac{z^{1-\gamma }}{\gamma -1}-np(t)\ln
\left( 1-p(t)z+...\right) \right) .
\end{multline*}%
Finally, for the factor containing $u$, replace $F$ by the substitution $%
\overline{F}(x)=\rho z/t$, i.e. $x=U(t/(z\rho ))$:%
\begin{eqnarray*}
\int_{U(1/(p(t)z))}^{\infty }e^{-ux/U(t)}dF(x) &\sim &\frac{z}{t}%
\int_{0}^{tp(t)}e^{-uU(t/(z\rho ))/U(t)}d\rho \\
&\sim &\frac{z}{t}\int_{0}^{tp(t)}\left( e^{-uU(t/(z\rho ))/U(t)}-1\right)
d\rho +p(t)z \\
&\sim &p(t)z\left( 1-\frac{1}{tp(t)}\int_{0}^{\infty }(1-e^{-u(z\rho
)^{-\gamma }})d\rho \right) \\
&\sim &p(t)z\left( 1-\frac{1}{tp(t)}\int_{0}^{\infty }\frac{%
(1-e^{-uz^{-\gamma }\eta })}{\eta ^{1+1/\gamma }}d\eta \right)
\end{eqnarray*}%
so that, as $t\rightarrow \infty $,%
\begin{multline*}
\left( \int_{U(1/(p(t)z))}^{\infty }e^{-ux/(tp(t)U(p^{-1}(t))}dF(x)\right)
^{np(t)} \\
\sim \exp \left( np(t)\ln z+np(t)\ln p(t)-\frac{n}{t}\int_{0}^{\infty }\frac{%
(1-e^{-uz^{-\gamma }\eta })}{\eta ^{1+1/\gamma }}d\eta \right).
\end{multline*}%

For the factor with the factorials, we have by Stirling's formula%
\begin{equation}  \label{stirl}
\frac{n!}{np(t)!(n(1-p(t))-1)!}\sim \frac{1}{\sqrt{2\pi }}\frac{n^{1/2}}{%
e^{(np(t)+1/2)\ln (p(t))}e^{(n(1-p(t))+1/2)\ln (1-p(t)}}.
\end{equation}
Equivalent for the integral in $z$:%
\begin{eqnarray*}
g(z) &=&\ln (z)-z \\
g^{\prime }(z) &=&\frac{1}{z}-1\qquad g^{\prime }(1)=0 \\
g^{\prime \prime }(z) &=&-\frac{1}{z^{2}}\qquad g^{\prime \prime }(1)=-1
\end{eqnarray*}%
By Laplace's method, we deduce that%
\begin{equation}  \label{lapl}
\int_{0}^{\infty }\exp \left( np(t)(\ln (z)-z)\right) dz\sim \frac{\sqrt{%
2\pi }}{n^{1/2}p^{1/2}(t)}\exp \left( -np(t)\right).
\end{equation}
Altogether%
\begin{eqnarray*}
&& \left.\Omega_{s}(u/U(t),v/U(p^{-1}(t)),w/(tp(t)U(p^{-1}(t))))\right\vert
_{N(t)=n} \\
&\sim &e^{-v}\exp \left( -\frac{n}{t}w\frac{1}{\gamma -1}+o\left( \frac{n}{t}%
\right) +np(t)\ln \left( 1-p(t)\right) +np(t)\ln p(t)\right) \\
&&\times \exp \left( -\frac{n}{t}\int_{0}^{\infty }\frac{(1-e^{-uz^{-\gamma
}\eta})}{\eta^{1+1/\gamma }}d\eta-np(t)\ln p(t)-np(t)-\frac{1}{2}\ln
(p(t))+\ln (p(t))\right) \\
&&\times \exp \left( -(np(t)+1/2)\ln (p(t))-(n(1-p(t))+1/2)\ln (1-p(t)\right)
\\
&\sim &\exp \left( -\frac{n}{t}\int_{0}^{\infty }\frac{(1-e^{-uz^{-\gamma
}\eta})}{\eta^{1+1/\gamma }}d\eta-v-\frac{n}{t}w\frac{1}{\gamma -1}\right).
\end{eqnarray*}%
\hfill $\Box$\newline

\textbf{Proof of Proposition \ref{Prop3}:} Again, we condition on $N(t)=n$:%
\begin{eqnarray*}
&&E_{N(t)=n}\left\{ \exp (-u\Lambda _{s}(t)/U(t)-vX_{n(1-p)}^{\ast }-w\Sigma
_{s}(t))/t\right\} \\
&=&\frac{n!}{np!(n(1-p)-1)!}\int_{0}^{\infty }\left( \int_{y}^{\infty
}e^{-ux/U(t)}dF(x)\right) ^{np}e^{-vy}\left(
\int_{0}^{y}e^{-wx/t}dF(x)\right) ^{n(1-p)-1}dF(y).
\end{eqnarray*}%
For the factor containing $w$, we have 
\begin{eqnarray*}
\int_{0}^{y}e^{-wx/t}dF(x) &=&P(X\leq y)-\frac{w}{t}%
\int_{0}^{y}xdF(x)+o(t^{-1}) \\
&=&P(X\leq y)\left( 1-\frac{w}{t}E\left\{ X|X\leq y\right\} +o(t^{-1})\right)
\\
&=&F(y)\left( 1-\frac{w}{t}E\left\{ X|X\leq y\right\} +o(t^{-1})\right) .
\end{eqnarray*}%
For the factor involving $u$ one can write 
\begin{eqnarray*}
\int_{y}^{\infty }e^{-ux/U(t)}dF(x) &=&-\left[ e^{-ux/U(t)}\overline{F}(x)%
\right] _{y}^{\infty }-\int_{y}^{\infty }e^{-ux/U(t)}\frac{x}{U(t)}\overline{%
F}(x)dx \\
&=&e^{-uy/U(t)}\overline{F}(y)-\int_{uy/U(t)}^{\infty }e^{-w}\overline{F}%
\left( \frac{w}{u}U(t)\right) dw \\
&=&\overline{F}(y)-\frac{1}{t}u^{\alpha }\int_{0}^{\infty }e^{-w}w^{-\alpha
}dw\left( 1+o(1)\right) \\
&=&\overline{F}(y)-\frac{1}{t}u^{\alpha }\Gamma (1-\alpha )\left(
1+o(1)\right) .
\end{eqnarray*}%
The ratio with factorials behaves, by Stirling's formula, as%
\begin{eqnarray*}
\frac{n!}{np!(n(1-p)-1)!} &\sim &\frac{\sqrt{2\pi }n^{n+1/2}e^{-n}}{\sqrt{%
2\pi }(np)^{np+1/2}e^{-np}\sqrt{2\pi }(n(1-p))^{n(1-p)+1/2}e^{-np}} \\
&\sim &\frac{1}{\sqrt{2\pi }}\frac{n^{1/2}}{p^{np+1/2}(1-p)^{n(1-p)+1/2}}.
\end{eqnarray*}%
For the integral in $y$ we have the equivalence%
\begin{eqnarray*}
&&\left( \int_{y}^{\infty }e^{-ux/t}dF(x)\right) ^{np}\left(
\int_{0}^{y}e^{-wx/t}dF(x)\right) ^{n(1-p)-1} \\
&=&\exp \left( n\left( p\ln (F(y)+(1-p)\ln (\overline{F}(y))\right) \right)
\exp \left( -w\frac{n}{t}pE\left\{ X|X\leq y\right\} -u\frac{n}{t}E\left\{
X|X>y\right\} \right) .
\end{eqnarray*}%
Let%
\begin{eqnarray*}
g(y) &=&p\ln (F(y))+(1-p)\ln (\overline{F}(y)) \\
g^{\prime }(y) &=&p\frac{f(y)}{F(y)}-(1-p)\frac{f(y)}{\overline{F}(y)} \\
F(y_{p}) &=&p,\qquad y_{p}=F^{-1}(p) \\
g(y_{p}) &=&p\ln (p)+(1-p)\ln (1-p) \\
g^{\prime \prime }(y) &=&p\frac{f^{\prime }(y)F(y)-f^{2}(y)}{F^{2}(y)}-(1-p)%
\frac{f^{\prime }(y)\overline{F}(y)+f^{2}(y)}{\overline{F}^{2}(y)} \\
g^{\prime \prime }(y_{p}) &=&f^{\prime }(y_{p})-\frac{f^{2}(y_{p})}{p}%
-f^{\prime }(y_{p})-\frac{f^{2}(y_{p})}{1-p}=-\frac{f^{2}(y_{p})}{p(1-p)}.
\end{eqnarray*}%
By Laplace's method, we deduce that%
\begin{equation*}
\exp \left( n\left( p\ln (F(y)+(1-p)\ln (\overline{F}(y))\right) \right)
\sim \frac{\sqrt{2\pi }\sqrt{p(1-p)}}{n^{1/2}}e^{n\left( p\ln (p)+(1-p)\ln
(1-p)\right) }
\end{equation*}%
and then%
\begin{eqnarray*}
&&\int_{0}^{\infty }\left( \int_{y}^{\infty }e^{-ux/t}dF(x)\right)
^{np}e^{-vy}\left( \int_{0}^{y}e^{-wx/t}dF(x)\right) ^{n(1-p)-1}dF(y) \\
&\sim &\frac{\sqrt{2\pi }\sqrt{p(1-p)}}{n^{1/2}}e^{n\left( p\ln (p)+(1-p)\ln
(1-p)\right) }\exp \left( -u\Theta E\left\{ X|X>y_{p}\right\}
-vy_{p}-w\Theta pE\left\{ X|X\leq y_{p}\right\} \right) .
\end{eqnarray*}%
Altogether%
\begin{equation*}
\Omega _{s}(u/t,v,w/t;t)\rightarrow e^{-vy_{p}}q_{0}\left( uE\left\{
X|X>y_{p}\right\} +wE\left\{ X|X\leq y_{p}\right\} \right) .
\end{equation*}%
\hfill $\Box $\newline

\textbf{Proof of Proposition \ref{Prop4}:} We first replace $F$ by the
substitution $\overline{F}(y)=z/t$, i.e. $y=U(t/z)$ \textbf{\ }\ 
\begin{multline*}
\Omega _{s}(u/U(t),v/U(t),w/t;t) =\sum_{n=0}^{s}p_{n}(t)\left(
\int_{0}^{\infty }e^{-ux/U(t)}dF(x)\right) ^{s} \\
+\frac{1}{s!}\int_{0}^{t}\left( t\int_{U(t/z)}^{\infty
}e^{-ux/U(t)}dF(x)\right) ^{s}e^{-vU(t/z)/U(t)}\frac{1}{t^{s+1}}%
Q_{t}^{(s+1)}\left( \int_{0}^{U(t/z)}e^{-wx/t}dF(x)\right) dz.
\end{multline*}%
Then we have%
\begin{eqnarray*}
\int_{0}^{U(t/z)}e^{-wx/t}dF(x) &=&F(U(t/z))-\frac{w}{t}%
\int_{0}^{U(t/z)}xdF(x)+o\left( \frac{1}{t}\right) \\
&=&1-\frac{z}{t}-\frac{w}{t}\mu +o\left( \frac{1}{t}\right)= 1-\frac{1}{t}%
\left( z+w\mu +o\left( 1\right) \right) .
\end{eqnarray*}%
Now use analogous arguments as in the proof of Proposition \ref{Prop1} and
note that 
\begin{eqnarray*}
E\left\{ \exp (-uS(t)/t)\right\} &=&Q_{t}\left( E\left\{ \exp
(-uX/t)\right\} \right) =Q_{t}\left( \int_{0}^{\infty }e^{-ux/t}dF(x)\right)
\\
&=&Q_{t}\left( 1-\frac{u}{t}\int_{0}^{\infty }xdF(x)+o\left( \frac{1}{t}%
\right) \right) \\
&\rightarrow &q_0\left( u\mu \right).
\end{eqnarray*}%
\hfill $\Box$\newline

\textbf{Proof of Corollary \ref{Prop4moment}:} By Proposition \ref{Prop4}

\begin{equation*}
E\left\{ \exp (-v\Xi _{s}-u\Sigma _{s})\right\} =\frac{1}{s!}%
\int_{0}^{\infty }z^{s}e^{-vz^{-\gamma }}q_{s+1}\left( z+u\mu \right) dz.
\end{equation*}%
Hence%
\begin{equation*}
\left. \frac{\partial }{\partial u}E\left\{ \exp (-v\Xi _{s}-u\Sigma
_{s})\right\} \right\vert _{u=0}=-\frac{\mu }{s!}\int_{0}^{\infty
}z^{s}e^{-vz^{-\gamma }}q_{s+2}\left( z\right) dz
\end{equation*}%
and therefore%
\begin{eqnarray*}
-\int_{0}^{\infty }\left. \frac{\partial }{\partial u}E\left\{ \exp
(-(u+v)\Xi _{s}-u\Sigma _{s})\right\} \right\vert _{u=0}dv &=&\frac{\mu }{s!}%
\int_{0}^{\infty }z^{s-\gamma }q_{s+2}\left( z\right) dz \\
&=&\mu \frac{\Gamma (s-\gamma +1)}{s!}E\left\{ \Theta ^{1+\gamma }\right\} .
\end{eqnarray*}
By Proposition \ref{Prop4} 
\begin{equation*}
E\left\{ \exp (-(u+v)\Xi _{0}-u\Sigma _{0})\right\} =\int_{0}^{\infty
}e^{-(u+v)z^{-\gamma }}q_{1}\left( z+u\mu \right) dz
\end{equation*}%
\begin{equation*}
\left. \frac{\partial }{\partial v}E\left\{ \exp (-(u+v)\Xi _{0}-u\Sigma
_{0})\right\} \right\vert _{v=0}=-\int_{0}^{\infty }z^{-\gamma
}e^{-uz^{-\gamma }}q_{1}\left( z+u\mu \right) dz
\end{equation*}%
\begin{eqnarray*}
&&-\int_{0}^{\infty }\left. \frac{\partial }{\partial v}E\left\{ \exp
(-(u+v)\Xi _{0}-u\Sigma _{0})\right\} \right\vert _{v=0}du \\
&=&\int_{0}^{\infty }\int_{0}^{\infty }z^{-\gamma }e^{-uz^{-\gamma
}}q_{1}\left( z+u\mu \right) du\,dz \\
&=&\int_{0}^{\infty }\left( -\left[ e^{-uz^{-\gamma }}q_{1}\left( z+u\mu
\right) \right] _{0}^{\infty }-\mu \int_{0}^{\infty }e^{-uz^{-\gamma
}}q_{2}\left( z+u\mu \right) du\right) dz \\
&=&\int_{0}^{\infty }q_{1}\left( z\right) dz-\mu \int_{0}^{\infty
}\int_{0}^{\infty }e^{-uz^{-\gamma }}q_{2}\left( z+u\mu \right) du\,dz \\
&=&1-\mu \int_{0}^{\infty }\int_{0}^{\infty }e^{-uz^{-\gamma }}q_{2}\left(
z+u\mu \right) du\,dz.
\end{eqnarray*}%
\hfill $\Box$\newline

\textbf{Proof of Proposition \ref{Prop5}:} Condition on $N(t)=n$ to see that 
\begin{eqnarray*}
&&E_{N(t)=n}\left\{ \exp (-u\Lambda
_{s}(t)/(tp(t)U(p^{-1}(t)))-vX_{n(1-p(t))}^{\ast }/U(p^{-1}(t)))-w\Sigma
_{s}(t))/t\right\} \\
&=&\frac{n!}{np(t)!(n(1-p(t)))!}\int_{0}^{\infty }\left( \int_{y}^{\infty
}e^{-ux/(tp(t)U(p^{-1}(t))}dF(x)\right) ^{np(t)}e^{-vy/U(p^{-1}(t))} \\
&&\times \left( \int_{0}^{y}e^{-wx/t}dF(x)\right) ^{n(1-p(t))-1}dF(y).
\end{eqnarray*}%
Now replace $F$ by the substitution $\overline{F}(y)=p(t)z$, i.e. $%
y=U(1/(p(t)z))$ \textbf{\ }\ 
\begin{eqnarray*}
&&\Omega _{s}(u/(tp(t)U(p^{-1}(t))),v/U(p^{-1}(t)),w/t;t) \\
&=&\frac{n!}{np(t)!(n(1-p(t)))!}\int_{0}^{\infty }\left(
\int_{U(1/(p(t)z))}^{\infty }e^{-ux/(tp(t)U(p^{-1}(t))}dF(x)\right) ^{np(t)}
\\
&&\times e^{-vU(1/(p(t)z))/U(p^{-1}(t))}\left(
\int_{0}^{U(1/(p(t)z))}e^{-wx/t}dF(x)\right) ^{n(1-p(t))-1}p(t)dz.
\end{eqnarray*}
Like before, 
\begin{equation*}
e^{-vU(1/(p(t)z))/U(p^{-1}(t))}\rightarrow e^{-vz^{-\gamma }}.
\end{equation*}
For the factor with $w$, one sees that 
\begin{eqnarray*}
\int_{0}^{U(1/(p(t)z))}e^{-wx/t}dF(x) &=&F(U(1/(p(t)z)))-\frac{w}{t}%
\int_{0}^{U(1/(p(t)z))}xdF(x)+o\left( \frac{1}{t}\right) ... \\
&=&1-p(t)z-\frac{w}{t}\mu +o\left( \frac{1}{t}\right) ...
\end{eqnarray*}%
and then%
\begin{eqnarray*}
\left( \int_{0}^{U(1/(p(t)z))}e^{-wx/t}dF(x)\right) ^{n(1-p(t))-1} &=&\exp
\left( n(1-p(t))\ln \left( 1-p(t)z-\frac{w}{t}\mu +o\left( \frac{1}{t}%
\right) \right) \right) \\
&=&\exp \left( -np(t)z-\frac{n}{t}w\mu +o\left( \frac{n}{t}\right) +np(t)\ln
\left( 1-p(t)z\right) \right).
\end{eqnarray*}

For the factor with $u$, replace $F$ by the substitution $\overline{F}%
(x)=p(t)\rho z$, i.e. $x=U(1/(p(t)z\rho ))$:%
\begin{equation*}
t\int_{U(t/z)}^{\infty }e^{-ux/U(t)}dF(x)=z\int_{0}^{1}e^{-uU(t/(z\rho
))/U(t)}d\rho \rightarrow z\int_{0}^{1}e^{-u(z\rho )^{-\gamma }}d\rho =\frac{%
z}{\gamma }\int_{1}^{\infty }\frac{e^{-uz^{-\gamma }\eta}}{\eta^{1+1/\gamma }%
}d\eta
\end{equation*}%
\begin{eqnarray*}
\int_{U(1/(p(t)z))}^{\infty }e^{-ux/(tp(t)U(p^{-1}(t))}dF(x)
&=&p(t)z\int_{0}^{1}e^{-uU(p^{-1}(t)/(z\rho ))/U(p^{-1}(t))tp(t)}d\rho \\
&=&p(t)z\left( 1-\frac{1}{tp(t)}\int_{0}^{1}u(z\rho )^{-\gamma }d\rho
+\right) \\
&=&p(t)z\left( 1-\frac{1}{1-\gamma }\frac{z^{-\gamma }}{tp(t)}%
u\int_{0}^{1}\rho ^{-\gamma }d\rho +\ldots\right)
\end{eqnarray*}%
and then%
\begin{eqnarray*}
&&\left( \int_{U(1/(p(t)z))}^{\infty }e^{-ux/(tp(t)U(p^{-1}(t))}dF(x)\right)
^{np(t)} \\
&=&\exp \left( np(t)\ln p(t)+np(t)\ln (z)-\frac{n}{t}\frac{u}{1-\gamma }%
z^{-\gamma }+\ldots\right).
\end{eqnarray*}
The ratio of factorials coincides with \eqref{stirl}. Also \eqref{lapl}
applies here. %For the integral in $z$: let%
% \begin{eqnarray*}
% g(z) &=&\ln (z)-z \\
% g^{\prime }(z) &=&\frac{1}{z}-1\qquad g^{\prime }(1)=0 \\
% g^{\prime \prime }(z) &=&-\frac{1}{z^{2}}\qquad g^{\prime \prime }(1)=-1
% \end{eqnarray*}%
% By Laplace's method%
% \begin{equation*}
% \int_{0}^{\infty }\exp \left( np(t)(\ln (z)-z)\right) dz\sim \frac{\sqrt{%
% 2\pi }}{n^{1/2}p^{1/2}(t)}\exp \left( -np(t)\right)
% \end{equation*}
Altogether%
\begin{eqnarray*}
&&\frac{n!}{np(t)!(n(1-p(t)))!}\int_{0}^{\infty }\left(
\int_{U(1/(p(t)z))}^{\infty }e^{-ux/(tp(t)U(p^{-1}(t))}dF(x)\right)
^{np(t)}e^{-vU(1/(p(t)z))/U(p^{-1}(t))} \\
&&\times \left( \int_{0}^{U(1/(p(t)z))}e^{-wx/t}dF(x)\right)
^{n(1-p(t))-1}p(t)dz \\
&\sim &e^{-v}\exp \left( -\frac{n}{t}w\mu +o\left( \frac{n}{t}\right)
+np(t)\ln \left( 1-p(t)\right) +np(t)\ln p(t)-\frac{n}{t}\frac{u}{1-\gamma }%
-np(t)-\frac{1}{2}\ln p(t)+\ln p(t)\right) \\
&&\times\exp \big( -(np(t)+1/2)\ln p(t)-(n(1-p(t))+1/2)\ln (1-p(t)\big) \\
&\sim &\exp \left( -\frac{n}{t}\frac{u}{1-\gamma }-v-\frac{n}{t}w\mu \right).
\end{eqnarray*}%
\hfill $\Box$\newline

\textbf{Proof of Proposition \ref{Prop6}:} Given $N(t)=n$%
\begin{eqnarray*}
&&E_{N(t)=n}\left\{ \exp (-u\Lambda _{s}(t)/t-vX_{n(1-p)}^{\ast }-w\Sigma
_{s}(t))/t\right\} \\
&=&\frac{n!}{np!(n(1-p)-1)!}\int_{0}^{\infty }\left( \int_{y}^{\infty
}e^{-ux/t}dF(x)\right) ^{np}e^{-vy}\left( \int_{0}^{y}e^{-wx/t}dF(x)\right)
^{n(1-p)-1}dF(y).
\end{eqnarray*}

The part involving $w$ coincides with the one in the proof of Proposition %
\ref{Prop3}. % : $e^{-vy}$
% 
% ii) Part with $w$:%
% \begin{eqnarray*}
% \int_{0}^{y}e^{-wx/t}dF(x) &=&P(X\leq y)-\frac{w}{t}%
% \int_{0}^{y}xdF(x)+o(t^{-1}) \\
% &=&P(X\leq y)\left( 1-\frac{w}{t}E\left\{ X|X\leq y\right\} +o(t^{-1})\right)
% \\
% &=&F(y)\left( 1-\frac{w}{t}E\left\{ X|X\leq y\right\} +o(t^{-1})\right)
% \end{eqnarray*}
For the factor involving $u$, we have%
\begin{eqnarray*}
\int_{y}^{\infty }e^{-ux/t}dF(x) &=&P(X>y)-\frac{u}{t}\int_{y}^{\infty
}xdF(x)+o(t^{-1}) \\
&=&\overline{F}(y)\left( 1-\frac{u}{t}E\left\{ X|X>y\right\}
+o(t^{-1})\right)
\end{eqnarray*}
The rest of the proof is completely analogous to the one for Proposition \ref%
{Prop3}.\hfill $\Box$\newline

\textbf{Proof of Proposition \ref{Prop7}:} Use the substitution $\overline{F}%
(y)=z/t$, i.e. $y=U(t/z)$:%
\begin{eqnarray*}
&&\Omega _{s}^{(\mu )}(u\Lambda _{s}(t)/U(t),vX_{N(t)-s}^{\ast
}/U(t),w\Sigma _{s}^{(\mu )}(t)/U(t);t) \\
&=&\sum_{n=0}^{s}p_{n}(t)\left( \int_{0}^{\infty }e^{-ux/U(t)}dF(x)\right)
^{s} \\
&&+\frac{1}{s!}\int_{0}^{t}\left( t\int_{U(t/z)}^{\infty
}e^{-ux/U(t)}dF(x)\right) ^{s}e^{-vU(t/z)/U(t)} \\
&&\times \frac{1}{t^{s+1}}Q_{t}^{(s+1)}\left( 1-\frac{1}{t}\left( t-te^{w\mu
/U(t)}\int_{0}^{U(t/z)}e^{-wx/U(t)}dF(x)\right) \right) dz
\end{eqnarray*}%
We then replace $F$ by the substitution $\overline{F}(x)=\rho z/t$, i.e. $%
x=U(t/(z\rho ))$:%
\begin{eqnarray*}
&&t-te^{w\mu /U(t)}\int_{0}^{U(t/z)}e^{-wx/U(t)}dF(x) \\
&=&t+te^{w\mu /U(t)}\int_{0}^{U(t/z)}\left( 1-\frac{wx}{U(t)}%
-e^{-wx/U(t)}\right) dF(x)-te^{w\mu /U(t)}F(U(t/z)) \\
&&+te^{w\mu /U(t)}\int_{0}^{U(t/z)}\frac{wx}{U(t)}\,dF(x) \\
&=&t\left( 1-e^{w\mu /U(t)}\left( 1-\frac{z}{t}\right) \right)
+z\int_{1}^{\infty }\left( 1-w\frac{U(t/(z\rho ))}{U(t)}-e^{-wU(t/(z\rho
))/U(t)}\right) d\rho \\
&&+\frac{tw}{U(t)}e^{w\mu /U(t)}\left( \mu -\frac{z^{1-\gamma }}{1-\gamma }%
\frac{U(t)}{t}(1+o(1))\right) \\
&=&t\left( 1-\left( 1+\frac{w\mu }{U(t)}+\frac{1}{2}\left( \frac{w\mu }{U(t)}%
\right) ^{2}+o\left( \frac{1}{U^{2}(t)}\right) \right) \left( 1-\frac{z}{t}%
\right) \right) \\
&&+z\int_{1}^{\infty }\left( 1-w\frac{U(t/(z\rho ))}{U(t)}-e^{-wU(t/(z\rho
))/U(t)}\right) d\rho +\frac{tw\mu }{U(t)}\left( 1+\frac{w\mu }{U(t)}%
+o\left( \frac{1}{U(t)}\right) \right) \left( 1+O(1/t)\right) \end{eqnarray*}
\begin{eqnarray*}&=&z+z\int_{1}^{\infty }\left( 1-w\frac{U(t/(z\rho ))}{U(t)}-e^{-wU(t/(z\rho
))/U(t)}\right) d\rho -\frac{z^{1-\gamma }}{1-\gamma }w+O\left( \frac{1}{U(t)%
}\right) +O\left( \frac{t}{U^{2}(t)}\right) \\
&\rightarrow &z\left( 1+\int_{1}^{\infty }(1-w(z\rho )^{-\gamma
}-e^{-w(z\rho )^{-\gamma }})d\rho -\frac{z^{-\gamma }}{1-\gamma }w\right) \\
&=&z\left( 1+\frac{1}{\gamma }\int_{0}^{1}\frac{1-wz^{-\gamma
}v-e^{-wz^{-\gamma }\eta}}{\eta^{1+1/\gamma }}d\eta-\frac{z^{-\gamma }}{%
1-\gamma }w\right).
\end{eqnarray*}%
Now use the same arguments as in the proof of Proposition \ref{Prop1}.\hfill 
$\Box$\newline

\textbf{Proof of Corollary \ref{Prop7moment}:} Note that

\begin{eqnarray*}
&&E\left\{ \exp (-(u+v)\Xi _{s}-u\Sigma _{s}^{(\mu )})\right\} \\
&=&\frac{1}{s!}\int_{0}^{\infty }z^{s}e^{-vz^{-\gamma }}e^{-uz^{-\gamma
}}q_{s+1}\left( z\left( 1+\frac{1}{\gamma }\int_{0}^{1}\frac{1-uz^{-\gamma
}v-e^{-uz^{-\gamma }\eta}}{\eta^{1+1/\gamma }}d\eta-\frac{z^{-\gamma }}{%
1-\gamma }u\right) \right) dz
\end{eqnarray*}%
\begin{eqnarray*}
&&\left. \frac{\partial }{\partial u}E\left\{ \exp (-(u+v)\Xi _{s}-u\Sigma
_{s})\right\} \right\vert _{u=0} \\
&=&-\frac{1}{s!}\int_{0}^{\infty }z^{s}z^{-\gamma }e^{-vz^{-\gamma
}}q_{s+1}\left( z\right) dz+\frac{1}{\left( 1-\gamma \right) s!}%
\int_{0}^{\infty }z^{1+s}z^{-\gamma }e^{-vz^{-\gamma }}q_{s+2}\left(
z\right) dz
\end{eqnarray*}%
\begin{equation*}
-\int_{0}^{\infty }\left. \frac{\partial }{\partial u}E\left\{ \exp
(-(u+v)\Xi _{s}-u\Sigma _{s}^{(\mu )})\right\} \right\vert _{u=0}dv=1+\frac{%
s+1}{ \gamma -1 }.
\end{equation*}%
\hfill $\Box$\newline

\textbf{Proof of Proposition \ref{Prop8}:} Use the substitution $\overline{F}%
(y)=z/t$, i.e. $y=U(t/z)$:%
\begin{eqnarray*}
&&\Omega _{s}^{(\mu )}(u\Lambda _{s}(t)/U(t),vX_{N(t)-s}^{\ast
}/U(t),w\Sigma _{s}^{(\mu )}/t^{1/2};t) \\
&=&\sum_{n=0}^{s}p_{n}(t)\left( \int_{0}^{\infty }e^{-ux/U(t)}dF(x)\right)
^{s} \\
&&+\frac{1}{s!}\int_{0}^{t}\left( t\int_{U(t/z)}^{\infty
}e^{-ux/U(t)}dF(x)\right) ^{s}e^{-vU(t/z)/U(t)} \\
&&\times \frac{1}{t^{s+1}}\,Q_{t}^{(s+1)}\left( 1-\frac{1}{t}\left(
t-te^{w\mu /t^{1/2}}\int_{0}^{U(t/z)}e^{-wx/t^{1/2}}dF(x)\right) \right) dz.
\end{eqnarray*}%
Then we have%
\begin{eqnarray*}
&&t-te^{w\mu /t^{1/2}}\int_{0}^{U(t/z)}e^{-wx/t^{1/2}}dF(x) \\
&=&t+te^{w\mu /t^{1/2}}\int_{0}^{U(t/z)}\left( 1-\frac{wx}{t^{1/2}}+\frac{1}{%
2}\frac{\left( wx\right) ^{2}}{t}-e^{-wx/t^{1/2}}\right) dF(x)-te^{w\mu
/t^{1/2}}F(U(t/z)) \\
&&+te^{w\mu /t^{1/2}}\int_{0}^{U(t/z)}\frac{wx}{t^{1/2}}dF(x)-\frac{1}{2}%
e^{w\mu /t^{1/2}}\int_{0}^{U(t/z)}\left( wx\right) ^{2}dF(x).
\end{eqnarray*}%
First note that%
\begin{eqnarray*}
&&te^{w\mu /t^{1/2}}\int_{0}^{U(t/z)}\left( 1-\frac{wx}{t^{1/2}}+\frac{1}{2}%
\frac{\left( wx\right) ^{2}}{t}-e^{-wx/t^{1/2}}\right) dF(x) \\
&=&ze^{w\mu /t^{1/2}}\int_{1}^{\infty }\left( 1-w\frac{U(t/(z\rho ))}{t^{1/2}%
}+\frac{1}{2}\left( w\frac{U(t/(z\rho ))}{t^{1/2}}\right)
^{2}-e^{-wU(t/(z\rho ))/t^{1/2}}\right) d\rho \rightarrow 0.
\end{eqnarray*}%
Secondly,%
\begin{eqnarray*}
&&t-te^{w\mu /t^{1/2}}F(U(t/z))+te^{w\mu /t^{1/2}}\int_{0}^{U(t/z)}\frac{wx}{%
t^{1/2}}dF(x)-\frac{1}{2}e^{w\mu /t^{1/2}}\int_{0}^{U(t/z)}\left( wx\right)
^{2}dF(x) \\
&=&t\left( 1-e^{w\mu /t^{1/2}}\left( 1-\frac{z}{t}\right) \right) +\frac{tw}{%
t^{1/2}}e^{w\mu /t^{1/2}}\left( \mu -U(t/z)\frac{z}{t}-\int_{U(t/z)}^{\infty
}\overline{F}(x)dx\right) \\
&&-\frac{1}{2}e^{w\mu /t^{1/2}}\left( E\left\{ X_{1}^{2}\right\}
-\int_{U(t/z)}^{\infty }x^{2}dF(x)\right) \\
&=&t\left( 1-e^{w\mu /t^{1/2}}\left( 1-\frac{z}{t}\right) \right) +\frac{tw}{%
t^{1/2}}e^{w\mu /t^{1/2}}\left( \mu -\frac{\alpha }{\alpha -1}U(t/z)\frac{z}{%
t}(1+o(1))\right) \\
&=&-\frac{1}{2}w^{2}e^{w\mu /t^{1/2}}\left( E\left\{ X_{1}^{2}\right\}
-\gamma \left( U(t/z)\right) ^{2}\frac{z}{t}(1+o(1)\right) \\
&=&t\left( 1-\left( 1+\frac{w\mu }{t^{1/2}}+\frac{1}{2}\left( \frac{w\mu }{%
t^{1/2}}\right) ^{2}+o\left( \frac{1}{t}\right) \right) \left( 1-\frac{z}{t}%
\right) \right) \\
&& +\,t^{1/2}\mu w\left( 1+\frac{w\mu }{t^{1/2}}+\frac{1}{2}\left( \frac{%
w\mu }{t^{1/2}}\right) ^{2}+o\left( \frac{1}{t}\right) \right) \left( 1-%
\frac{\alpha }{\mu (\alpha -1)}U(t/z)\frac{z}{t}(1+o(1)\right) \\
&&-\frac{1}{2}w^{2}\left( 1+\frac{w\mu }{t^{1/2}}+\frac{1}{2}\left( \frac{%
w\mu }{t^{1/2}}\right) ^{2}+o\left( \frac{1}{t}\right) \right) \left(
E\left\{ X_{1}^{2}\right\} -\gamma \left( U(t/z)\right) ^{2}\frac{z}{t}%
(1+o(1)\right)
\end{eqnarray*}%
and it follows that%
\begin{eqnarray*}
&&t-te^{w\mu /t^{1/2}}F(U(t/z))+te^{w\mu /t^{1/2}}\int_{0}^{U(t/z)}\frac{wx}{%
t^{1/2}}dF(x)-\frac{1}{2}e^{w\mu /t^{1/2}}\int_{0}^{U(t/z)}\left( wx\right)
^{2}dF(x) \\
&=&z-t^{1/2}\mu w-\frac{1}{2}\left( w\mu \right) ^{2}+o(1) +t^{1/2}\mu
w+\left( w\mu \right) ^{2}+O\left( \frac{U(t)}{t^{1/2}}\right) -\frac{1}{2}%
w^{2}E\left\{ X_{1}^{2}\right\} +O\left( \left( \frac{U(t)}{t^{1/2}}\right)
^{2}\right) \\
&\rightarrow &z-\frac{1}{2}w^{2}\sigma ^{2}
\end{eqnarray*}%
which completes the proof. Note that%
\begin{eqnarray*}
E\left\{ \exp (-u(S(t)-N(t)\mu )/t^{1/2})\right\} &=&Q_{t}\left( E\left\{
\exp (-u(X-\mu )/t^{1/2})\right\} \right) \\
&=&Q_{t}\left( \int_{0}^{\infty }e^{-u(x-\mu )/t^{1/2}}dF(x)\right) \\
&=&Q_{t}\left( 1+\frac{u^{2}}{2t}\sigma ^{2}+o\left( \frac{1}{t}\right)
\right) . \\
&\rightarrow &q_0\left( -\frac{u^{2}}{2}\sigma ^{2}\right) =E\left\{
e^{u^{2}\sigma ^{2}\Theta /2}\right\}
\end{eqnarray*}%
\hfill $\Box$\newline

\section{Conclusion}

\label{secc}

In this paper we provided a fairly general collection of results on the
joint asymptotic Laplace transforms of the normalized sums of smallest and
largest among regularly varying claims, when the length of the considered
time interval tends to infinity. This extends several classical results in
the field. The appropriate scaling of the different quantities is essential.
We showed to what extent the type of the near mixed Poisson process counting
the number of claim instances influences the limit results, and also
identified quantities for which this influence is asymptotically negligible.
We further related the dominance of the maximum term in such a random sum to
another quantity that exhibits the effect of the tail index on the aggregate
claim rather explicitly, namely the ratio of sum of squares of the claims
over the sum of the claims squared. The results allow to further quantify
the effect of large claims on the total claim amount in an insurance
portfolio, and could hence be helpful in the design of appropriate
reinsurance programs when facing heavy-tailed claims with regularly varying
tail. Particular emphasis is given to the case when the tail index exceeds
1, which corresponds to infinite-mean claims, a situation that is
particularly relevant for catastrophe modelling.\\

% 
% \section{Conclusion}
% XXX.

%%%%%%%%%%%%%%%%%%%%%%%%%%%%%%%%%%%%%%%%%%

%%%%%%%%%%%%%%%%%%%%%%%%%%%%%%%%%%%%%%%%%%

%%%%%%%%%%%%%%%%%%%%%%%%%%%%%%%%%%%%%%%%%%

%%%%%%%%%%%%%%%%%%%%%%%%%%%%%%%%%%%%%%%%%%

%%%%%%%%%%%%%%%%%%%%%%%%%%%%%%%%%%%%%%%%%%

%%%%%%%%%%%%%%%%%%%%%%%%%%%%%%%%%%%%%%%%%%
% 
% \conflictofinterests{Conflicts of Interest}
% 
% State any potential conflicts of interest here or "The authors declare no conflicts of interest". 

%=================================================================
% References: Variant A
%=================================================================
% Back Matter (References and Notes)
%----------------------------------------------------------
% Style and layout of the references
\makeatletter
\renewcommand\@biblabel[1]{#1. } \makeatother

%=================================================================
% References:  Variant B
%=================================================================
% Use the following option to include external BibTeX files:
%\bibliography{lite}
%\bibliographystyle{mdpi}
\end{document}